\documentclass[reqno]{amsart}
\usepackage{amssymb}
\usepackage{amsmath}
\usepackage[latin1]{inputenc}
\usepackage{fullpage}
\makeatletter
\@addtoreset{equation}{section}
\makeatother

\renewcommand\thefigure{\thesection.\@arabic\c@figure}
\renewcommand\thetable{\thesection.\@arabic\c@table}

\newtheorem{theorem}{Theorem}[section]
\newtheorem{lemma}[theorem]{Lemma}
\newtheorem{proposition}[theorem]{Proposition}

\newcommand{\mf}[1]{{\mathfrak #1}}
\newcommand{\mc}[1]{{\mathcal #1}}
\newcommand{\bs}[1]{{\boldsymbol #1}}
\newcommand{\bb}[1]{{\mathbb #1}}
\begin{document}
\title[Hydrodynamic limit for a stochastic lattice gas model]{Hydrodynamic limit for a boundary driven stochastic lattice gas model with many conserved quantities}
\author{Alexandre B. Simas}
\address{\noindent IMPA, Estrada Dona Castorina 110, CEP 22460 Rio de
  Janeiro, Brasil. 
\newline e-mail: \rm
  \texttt{alesimas@impa.br} }

\thanks{Research supported by CNPq.}
\begin{abstract}
  We prove the hydrodynamic limit for a particle system in which particles may have different velocities. We assume that we have two infinite reservoirs of particles at the boundary: this is the so-called boundary driven process. The dynamics we considered consists of a weakly asymmetric simple exclusion process with collision among particles having different velocities. 
\end{abstract}
\maketitle
\section{Introduction}
Interacting particle systems have been the subject of intense studies during the
last 30 years due to the fact that, in one hand, they present many of the collective
features that are found in real physical systems, and, in the other hand they are, up
to some extent, mathematically tractable. Their study has led in many cases to a
more detailed understanding of the microscopic mechanisms behind those collective
phenomena. We refer to \cite{KL} for further references, and to \cite{BDGJL} for recent results.

Since their introduction by Spitzer \cite{S}, the simple exclusion process and the
zero-range process have been among the most studied interacting particles systems,
and they have served as a test field for new mathematical and physical ideas.

In the last years there has been considerable progress in understanding stationary non equilibrium states: reversible systems in contact
with different reservoirs at the boundary imposing a gradient on the conserved quantities of the system. In these systems there is a flow of matter
through the system and the dynamics is not reversible. The main difference with respect to equilibrium (reversible)
states is the following. In equilibrium, the invariant measure, which determines the
thermodynamic properties, is given for free by the Gibbs distribution specified by
the Hamiltonian. On the contrary, in non equilibrium states the construction of
the stationary state requires the solution of a dynamical problem. One of the most striking typical property of these systems is the presence of long-range correlations. For the symmetric simple exclusion this was already shown in a pioneering paper by Spohn \cite{S2}. 
We refer to \cite{BDGJL1,D} for two recent reviews on this topic. 

The hydrodynamic behavior of the one-dimensional boundary driven exclusion process was studied by \cite{ELS1}, \cite{ELS2} and \cite{KLO}. Also,
Landim, Olla and Volchan \cite{LOV} considered the behavior of a tagged particle in a one-dimensional nearest-neighbor symmetric exclusion process under
the action of an external constant, and made connections between the behavior of a tagged particle in this situation and a process with infinite reservoirs.

We consider a stationary non-equilibrium state, whose non-equilibrium is due to external
fields and/or different chemical potentials at the boundaries, in which there is a flow of
physical quantities, such as heat, electric charge, or chemical substances, across the
system. The hydrodynamic behavior for this kind of processes in any dimension has been
solved by \cite{ELS1,ELS2}. Nevertheless,
they have solved this problem only for the case where the unique thermodynamic observable quantity is 
the empirical density. 

Our goal is to extend their results to the situation when there are several thermodynamic variables:
density and momentum. It is not always clear that a closed macroscopic dynamical
description is possible. However, we show that the system can be described by a hydrodynamic equation:
fix a macroscopic time interval
$[0, T ]$, and consider the dynamical behavior of the empirical density and momentum over such an
interval. The law of large numbers for the empirical density and momentum is then called hydrodynamic
limit and, in the context of the diffusive scaling limit here considered, is given
by a system of parabolic evolution equations which is called hydrodynamic equation. Once the hydrodynamic limit for this model is rigorously established, a reasonable
goal is to find an explicit connection between the thermodynamic potentials and
the dynamical macroscopic properties like transport coefficients. The study of large
deviations provides such a connection. The dynamical large deviation for boundary driven exclusion processes in any dimension with one conserved quantity has been recently proved in \cite{flm}. 

The dynamical large deviations for the model with many conserved quantities is studied at \cite{fsv}, and the hydrodynamic limit obtained in this article is important for such large deviations.

The model which we will study can be informally described as follows: fix a velocity $v$, an integer $N\geq 1$, and boundary densities $0 < \alpha_v(\cdot) < 1$ and $0 < \beta_v(\cdot) < 1$; at any given time, each site of the set $\{1,\ldots,N-1\}\times\{0,\ldots,N-1\}^{d-1}$ is either empty or occupied by one particle at velocity $v$. In the bulk, each particle attempts to jump at any of its neighbors at the same velocity, with a weakly asymmetric rate. To respect the exclusion rule, the particle jumps only if the target
site at the same velocity $v$ is empty; otherwise nothing happens. At the boundary, sites with first coordinates given by $1$ or $N -1$ have particles
being created or removed in such a way that the local densities are $\alpha_v(\tilde{x})$ and $\beta_v(\tilde{x})$: at rate $\alpha_v(\tilde{x}/N)$ a particle is
created at $\{1\}\times\{\tilde{x}\}$ if the site is empty, and at rate $1-\alpha_v(\tilde{x})$ the particle at $\{1\}\times\{\tilde{x}\}$
is removed if the site is occupied, and at rate $\beta_v(\tilde{x})$ a particle is
created at $\{N-1\}\times\{\tilde{x}\}$ if the site is empty, and at rate $1-\beta_v(\tilde{x})$ the particle at $\{N-1\}\times\{\tilde{x}\}$
is removed if the site is occupied. Superposed to this dynamics, there is a collision process which exchange
velocities of particles in the same site in a way that momentum is conserved.

Similar models have been studied by \cite{BL,emy,qy}. In fact, the model we consider here is based on the model of Esposito et al. \cite{emy} which was used to derive the Navier-Stokes equation. It is also noteworthy that the derivation of hydrodynamic limits and macroscopic fluctuation theory for a system with two conserved quantities have been studied in \cite{cb}.

Under diffusive time scaling, assuming local equilibrium, it is not difficult to
show that the evolution of the thermodynamic quantities is described by the parabolic system of 
equations
\begin{equation}\label{eq1}
\partial_t (\rho,{\boldsymbol p}) + \sum_{v\in \mathcal{ V}} \tilde{v}\left[v\cdot \nabla F(\rho,{\boldsymbol p})\right] = \frac{1}{2}\Delta (\rho,{\boldsymbol p}),\\
\end{equation}
where $\tilde{v} = (1,v_1,\ldots,v_d)$, $\rho$ stands for the density and ${\boldsymbol p} = (p_1, \ldots , p_d)$ for the momentum. $F$ is a
thermodynamical quantity determined by the ergodic properties of the dynamics.

Therefore, the purpose of this article is to define an interacting particle system whose
macroscopic density profile evolves according to the partial differential equation given by (\ref{eq1}) with initial condition
$$(\rho,{\boldsymbol p})(0,\cdot) =  (\rho_0,{\boldsymbol p}_0)(\cdot) \hbox{~and~} (\rho,{\boldsymbol p})(t,x) = (\rho,{\boldsymbol p})_b(x), x\in \partial D,$$
with $D$ being a suitable domain, and the equality on the boundary being on the trace sense.

This equation derives from the underlying stochastic dynamics through an appropriate
scaling limit in which the microscopic time and space coordinates are
rescaled diffusively. The hydrodynamic equation (\ref{eq1}) thus represents the law of
large numbers for the empirical density and momentum of the stochastic lattice gas. The convergence
has to be understood in probability with respect to the law of the stochastic
lattice gas. Finally, the initial condition for (\ref{eq1}) depends on the initial distribution
of particles. Of course many microscopic configurations give rise to the same initial
condition $(\rho_0,{\boldsymbol p}_0)(\cdot)$.

The article is organized as follows: in Section \ref{sec2} we establish the notation and
state the main results of the article; in Section \ref{sec4}, we prove the hydrodynamic limit
for the particle system we are interested in; the proof of a Replacement Lemma needed
for the hydrodynamic limit is postponed to Section \ref{sec5}; in Section
\ref{sec6} we prove the uniqueness of weak solutions of the hydrodynamic equations also needed for the hydrodynamic limits.
\section{Notation and results}\label{sec2}
Let $\mathbb{T}_N^d = \{0,\ldots,N-1\}^d = (\mathbb{Z}/N\mathbb{Z})^d$, the $d$-dimensional discrete torus, and
let $D_N^d = S_N\times \mathbb{T}_N^{d-1}$, with $S_N = \{1,\ldots,N-1\}$. Further, let also $\mathcal{V}\subset \mathbb{R}^d$ be a finite set of velocities $v = (v_1,\ldots,v_d)$. Assume that $\mathcal{V}$ is invariant under reflexions and permutations of the coordinates:
$$(v_1,\ldots,v_{i-1},-v_i,v_{i+1},\ldots,v_d)\hbox{~and~}(v_{\sigma(1)},\ldots,v_{\sigma(d)})$$
belong to $\mathcal{V}$ for all $1\leq i\leq d$, and all permutations $\sigma$ of $\{1,\ldots,d\}$, provided $(v_1,\ldots,v_d)$ belongs to $\mathcal{ V}$. Finally, denote the $d$-dimensional torus by $\mathbb{T}^d = [0,1)^d = (\mathbb{R}/\mathbb{Z})^d$. 

On each site of $D_N^d$, at most one particle for each velocity is allowed. We denote: the number of particles with velocity $v$ at $x$, $v\in \mathcal{ V}$, $x\in D_N^d$, by $\eta(x,v)\in \{0,1\}$;
the number of particles in each velocity $v$ at a site $x$ by $\eta_x = \{\eta(x,v); v\in\mathcal{ V}\}$; and a configuration by $\eta = \{\eta_x; x\in D_N^d\}$. The set of particle configurations is $X_N = \left( \{0,1\}^\mathcal{V}\right)^{D_N^d}$.

On the interior of the domain, the dynamics consists of two parts: (i) each particle of the system evolves according to a nearest neighbor weakly asymmetric random walk with exclusion among particles of the same velocity, and (ii) binary collision between particles of different velocities. Let $p(x,v)$ be an irreducible probability transition function of finite range, and mean velocity $v$:
$$\sum_x x p(x,v) = v.$$
The jump law and the waiting times are chosen so that the jump rate from site $x$ to site $x+y$ for a particle with velocity $v$ is
$$P_N(y,v) = \frac{1}{2}\sum_{j=1}^d (\delta_{y,e_j} + \delta_{y,-e_j}) + \frac{1}{N}p(y,v),$$
where $\delta_{x,y}$ stands for the Kronecker delta, which equals one if $x=y$ and 0 otherwise, and $\{e_1,\ldots,e_d\}$ is the canonical basis in $\mathbb{R}^d$.

\subsection{The boundary driven exclusion process}
Our main interest is to examine the stochastic lattice gas model given by the
generator ${\mathcal{L}}_N$ which is
the superposition of the boundary dynamics with the collision and exclusion:
\begin{equation}\label{geradorbd}
{\mathcal{ L}}_N = N^2\{\mathcal{L}_N^b + \mathcal{ L}_N^c +\mathcal{ L}_N^{ex}\},
\end{equation}
where $\mathcal{L}_N^b$ stands for the generator which models the part of the dynamics at which a particle at the boundary can enter or leave the system, $\mathcal{ L}_N^c$ stands for the generator which models the collision part of the dynamics and lastly, $\mathcal{ L}_N^{ex}$ models the exclusion part of the dynamics. Let $f$ be a local function on $X_N$.
The generator of the exclusion part of the dynamics, $\mathcal{L}_N^{ex}$, is given by
$$(\mathcal{ L}_N^{ex} f)(\eta) = \sum_{v\in\mathcal{ V}} \sum_{x,x+z\in D_N^d} \eta(x,v)[1-\eta(z,v)]P_N(z-x,v)\left[f(\eta^{x,z,v})-f(\eta) \right],$$
where
$$\eta^{x,y,v}(z,w) = \left\{
\begin{array}{cc}
\eta(y,v)&\text{if $w=v$ and $z=x$},\\
\eta(x,v)&\text{if $w=v$ and $z=y$},\\
\eta(z,w)&\text{otherwise}.
\end{array}
\right.$$

We will often use the decomposition 
$$\mathcal{ L}_N^{ex} = \mathcal{ L}_N^{ex,1} + \mathcal{ L}_N^{ex,2},$$
where
$$(\mathcal{ L}_N^{ex,1} f)(\eta) = \frac{1}{2}\sum_{v\in\mathcal{ V}} \sum_{\substack{x,x+z\in D_N^d\\|z-x|=1}} \eta(x,v)[1-\eta(z,v)]\left[f(\eta^{x,z,v})-f(\eta) \right],$$
and
$$(\mathcal{ L}_N^{ex,2} f)(\eta) = \frac{1}{N}\sum_{v\in\mathcal{ V}} \sum_{x,x+z\in D_N^d} \eta(x,v)[1-\eta(z,v)]p(z-x,v)\left[f(\eta^{x,z,v})-f(\eta) \right].$$

The generator of the collision part of the dynamics, $\mathcal{ L}_N^c$, is given by
$$(\mathcal{ L}_N^c f)(\eta) = \sum_{y\in D_N^d}\sum_{q\in\mathcal{ Q}} p(y,q,\eta)\left[ f(\eta^{y,q}) - f(\eta)\right],$$
where $\mathcal{Q}$ is the set of all collisions which preserve momentum:
$$\mathcal{Q} = \{q=(v,w,v',w') \in \mathcal{ V}^4 : v+w=v'+w'\},$$
the rate $p(y,q,\eta)$ is given by
$$p(y,q,\eta) = \eta(y,v)\eta(y,w)[1-\eta(y,v')][1-\eta(y,w')],$$
and for $q = (v_0,v_1,v_2,v_3)$, the configuration $\eta^{y,q}$ after the collision is defined as
$$\eta^{y,q}(z,u) = \left\{
\begin{array}{cc}
\eta(y,v_{j+2})&\text{if $z=y$ and $u=v_j$ for some $0\leq j\leq 3$},\\
\eta(z,u)&\text{otherwise,}
\end{array}
\right.$$
where the index of $v_{j+2}$ should be taken modulo 4. 

Particles of velocities $v$ and $w$ at the same site collide at rate one and produce two particles of velocities $v'$ and $w'$ at that site.

Finally, the generator of the boundary part of the dynamics is given by
\begin{eqnarray*}
(\mathcal{ L}_N^b f)(\eta) &=&\!\!\! \sum_{\substack{x\in D_N^d\\x_1 = 1}} \sum_{v\in\mathcal{ V}} [ \alpha_v(\tilde{x}/N)[1-\eta(x,v)] + (1-\alpha_v(\tilde{x}/N))\eta(x,v)][f(\sigma^{x,v}\eta)-f(\eta)]\\
&+&\!\!\! \sum_{\substack{x\in D_N^d\\x_1 = N-1}} \sum_{v\in\mathcal{ V}} [ \beta_v(\tilde{x}/N)[1-\eta(x,v)] + (1-\beta_v(\tilde{x}/N))\eta(x,v)][f(\sigma^{x,v}\eta)-f(\eta)],
\end{eqnarray*} 
where $\tilde{x} = (x_2,\ldots,x_d)$, 
$$\sigma^{x,v}\eta(y,w) = \left\{\begin{array}{cc}
1-\eta(x,w),& \hbox{if ~} w = v\hbox{~and~} y=x,\\
\eta(y,w),& \hbox{otherwise.}
\end{array}
\right.,$$
and for every $v\in\mathcal{V}$, $\alpha_v,\beta_v\in C^2(\mathbb{T}^{d-1})$. We also assume that, for every $v\in\mathcal{V}$, $\alpha_v$ and $\beta_v$ have images belonging to some compact subset of $(0,1)$. The functions $\alpha_v$ and $\beta_v$, which
affect the birth and death rates at the two boundaries, represent the densities of the reservoirs.

Note that time has been speeded up diffusively in \eqref{geradorbd}. Let $\{\eta(t), t\geq 0\}$ be the Markov process with generator $\mathcal{ L}_N$ and denote by $\{ S_t^N, t\geq 0\}$ the semigroup associated to $\mathcal{ L}_N$.

Let $D(\mathbb{R}_+,X_N)$ be the set of right continuous functions with left limits taking values on $X_N$. For a probability measure $\mu$ on $X_N$, denote by $\mathbb{ P}_\mu$ the measure on the path space $D(\mathbb{ R}_{+}, X_N)$ induced by $\{\eta(t): t\geq 0\}$ and the initial measure $\mu$. Expectation with respect to $\mathbb{ P}_\mu$ is denoted by $\mathbb{ E}_\mu$.

\subsection{Mass and momentum}\label{massmomentum}
For each configuration $\xi \in\{0,1\}^\mathcal{ V}$, denote by $I_0(\xi)$ the mass of $\xi$ and by $I_k(\xi)$, $k=1,\ldots,d,$ the momentum of $\xi$:
$$I_0(\xi) = \sum_{v\in\mathcal{ V}} \xi(v),\quad I_k(\xi)=\sum_{v\in\mathcal{ V}} v_k \xi(v).$$
Set ${\boldsymbol I}(\xi) := (I_0(\xi),\ldots,I_d(\xi))$. Assume that the set of velocities is chosen in such a way that the unique quantities
conserved by the random walk dynamics described above are mass and momentum: $\sum_{x\in D_N^d} {\boldsymbol I}(\eta_x)$. Two examples of sets of velocities satisfying these conditions can be found at \cite{emy}.

For each chemical potential ${\boldsymbol \lambda} = (\lambda_0,\ldots,\lambda_d) \in \mathbb{ R}^{d+1}$, denote by $m_\mathbb{ \lambda}$ the measure on $\{0,1\}^\mathcal{ V}$ given by
\begin{equation}\label{mprod}
m_\mathbb{ \lambda} (\xi) = \frac{1}{Z({\boldsymbol \lambda})} \exp\left\{\mathbb{ \lambda}\cdot {\boldsymbol I}(\xi)\right\},
\end{equation}
where $Z({\boldsymbol \lambda})$ is a normalizing constant. Note that $m_{\boldsymbol \lambda}$ is a product measure on $\{0,1\}^\mathcal{ V}$, i.e., that the variables $\{\xi(v): v\in\mathcal{ V}\}$ are independent under $m_{\boldsymbol \lambda}$.

Denote by $\mu_{\boldsymbol \lambda}^N$ the product measure on $X_N$, with marginals given by 
\begin{equation*}
\mu_{\boldsymbol \lambda}^N \{\eta: \eta(x,\cdot) = \xi\} = m_{\boldsymbol \lambda}(\xi),
\end{equation*}
for each $\xi$ in $\{0,1\}^\mathcal{ V}$ and $x\in D_N^d$. Note that $\{\eta(x,v): x\in D_N^d, v\in\mathcal{ V}\}$ are independent variables under $\mu_{\boldsymbol \lambda}^N$,
and that the measure $\mu_{\boldsymbol \lambda}^N$ is invariant for the exclusion process with periodic boundary condition.

The expectation under $\mu_{\boldsymbol \lambda}^N$ of the mass and momentum are given by
\begin{eqnarray*}
\rho({\boldsymbol \lambda})&:=& E_{\mu_{\boldsymbol \lambda}^N} \left[ I_0(\eta_x)\right] = \sum_{v\in\mathcal{ V}} \theta_v({\boldsymbol \lambda}),\\
p_k({\boldsymbol \lambda}) &:=& E_{\mu_{\boldsymbol \lambda}^N} \left[ I_k(\eta_x)\right] = \sum_{v\in\mathcal{ V}} v_k \theta_v({\boldsymbol \lambda}).
\end{eqnarray*}
In this formula $\theta_v({\boldsymbol \lambda})$ denotes the expected value of the density of particles with velocity $v$ under $m_{\boldsymbol \lambda}$:
$$\theta_v({\boldsymbol \lambda}):= E_{m_{\boldsymbol \lambda}} \left[\xi(v)\right] = \frac{\exp\left\{\lambda_0 + \sum_{k=1}^d \lambda_k v_k\right\}}{1+ \exp\left\{\lambda_0 + \sum_{k=1}^d \lambda_k v_k\right\}}.$$

Denote by $(\rho,{\boldsymbol p})({\boldsymbol \lambda}):= (\rho({\boldsymbol \lambda}),p_1({\boldsymbol \lambda}),\ldots, p_d({\boldsymbol \lambda}))$ the map that associates the chemical potential to the vector of density and momentum. It is possible to prove that $(\rho,{\boldsymbol p})$ is a diffeomorphism onto ${\mf U}\subset \mathbb{ R}^{d+1}$, the interior of the convex envelope of $\left\{ {\boldsymbol I}(\xi), \xi\in\{0,1\}^\mathcal{ V}\right\}$. Denote by $\Lambda = (\Lambda_0,\ldots,\Lambda_d): {\mf U}\to\mathbb{ R}^{d+1}$ the inverse of $(\rho, {\boldsymbol p})$. This correspondence allows one to parameterize the invariant states by the density and momentum: for each $(\rho,{\boldsymbol p})$ in ${\mf U}$ we have a product measure $\nu_{\rho,{\boldsymbol p}}^N = \mu_{\Lambda(\rho,{\boldsymbol p})}^N$ on $X_N$.

\subsection{Hydrodynamic limit for the boundary driven exclusion process}
Let $D^d = [0,1]\times \bb T^{d-1}$. Fix $\rho_0:D^d\to\mathbb{R}_+$ and ${\boldsymbol p}_0:D^d\to\mathbb{R}^d$, where ${\boldsymbol p}_0=(p_{0,1},\ldots,p_{0,d})$. We say that a sequence of probability measures $(\mu_N)_N$ on $X_N$ is associated to the density profile $\rho_0$ and momentum profile ${\boldsymbol p}_0$, if, for every continuous function $G:D^d\to \mathbb{ R}$ and for every $\delta >0$,
$$\lim_{N\to\infty} \mu^{N} \left[ \eta: \left| \frac{1}{N^d} \sum_{x\in D^d_N} G\left(\frac{x}{N}\right) I_0(\eta_x) - \int_{D^d} G(u) \rho_0(u) du \right| > \delta\right] = 0,$$
and for every $1\leq k\leq d$
$$\lim_{N\to\infty} \mu^{N} \left[ \eta: \left| \frac{1}{N^d} \sum_{x\in D^d_N} G\left(\frac{x}{N}\right) I_k(\eta_x) - \int_{D^d} G(u) p_{0,k} (u) du\right| > \delta\right] =0.$$

Fix $T>0$ and let $(B,\|\cdot\|_B)$ be a Banach space. We denote by $L^2([0,T],B)$ the Banach space of measurable functions $U:[0,T]\to B$ for which
$$\|U\|_{L^2([0,T],B)}^2 = \int_0^T \|U_t\|_B^2 dt <\infty.$$
Moreover, we denote by $H^1(D^d)$ the Sobolev space of measurable functions in $L^2(D^d)$ that have generalized derivatives in $L^2(D^d)$.

For $x = (x_1,\tilde{x})\in \{0,1\}\times\mathbb{T}^{d-1}$, let
\begin{equation}\label{funcaod}
d(x) =\left\{
\begin{array}{cc}
a(\tilde{x}) = \sum_{v\in\mathcal{V}} (\alpha_v(\tilde{x}),v_1\alpha_v(\tilde{x}),\ldots,v_d\alpha_v(\tilde{x})),&\hbox{if~~}x_1=0,\\[10pt] 
b(\tilde{x}) = \sum_{v\in\mathcal{V}} (\beta_v(\tilde{x}),v_1\beta_v(\tilde{x}),\ldots,v_d\beta_v(\tilde{x})),&\hbox{if~~}x_1=1.
\end{array}
\right.
\end{equation}

Fix a bounded density profile $\rho_0:D^d \to \mathbb{ R}_{+}$, and a bounded momentum
profile ${\boldsymbol p}_0: D^d \to \mathbb{ R}^d$.
A bounded function $({\rho},{{\boldsymbol p}}): [0,T]\times D^d \to \mathbb{ R}_{+}\times \mathbb{ R}^d$ is a weak solution
of the system of parabolic partial differential equations
\begin{equation}\label{problema2}
\left\{ 
\begin{array}{c}
\partial_t (\rho,{\boldsymbol p}) + \sum_{v\in \mathcal{ V}} \tilde{v}\left[v\cdot \nabla\chi(\theta_v(\Lambda(\rho,{\boldsymbol p})))\right] = \frac{1}{2}\Delta (\rho,{\boldsymbol p}),\\[10pt]
(\rho,{\boldsymbol p})(0,\cdot) =  (\rho_0,{\boldsymbol p}_0)(\cdot) \hbox{~and~} (\rho,{\boldsymbol p})(t,x) = d(x), x\in \{0,1\}\times\mathbb{T}^{d-1},
\end{array}
\right.
\end{equation}
if for every vector valued function $H:[0,T]\times D^d\to\mathbb{ R}^{d+1}$ of class $C^{1,2}\left([0,T]\times D^d\right)$
vanishing at the boundary, we have
$$
\int_{D^d} H(T,u)\cdot ({\rho},{{\boldsymbol p}})(T,u)du - \int_{D^d} H(0,u)\cdot(\rho_0, {\boldsymbol p}_0) (u) du
$$
$$
=\int_0^T dt \int_{D^d} du \left\{ ({\rho},{{\boldsymbol p}})(t,u)\cdot\partial_t H(t,u) + \frac{1}{2} ({\rho},{{\boldsymbol p}})(t,u) \cdot \sum_{1\leq i\leq d} \partial_{u_i}^2 H(t,u)\right\}
$$
$$
-\int_0^T dt \int_{\{1\}\times\mathbb{ T}^{d-1}} dS\,\, b(\tilde{u})\cdot \partial_{u_1} H(t,u) + \int_0^T dt \int_{\{0\}\times\mathbb{ T}^{d-1}} dS\,\, a(\tilde{u})\cdot \partial_{u_1} H(t,u)
$$
$$
-\int_0^T dt \int_{D^d} du\,\,  \sum_{v\in\mathcal{ V}} \tilde{v}\cdot \chi(\theta_v(\Lambda({\rho},{{\boldsymbol p}}))) \sum_{1\leq i \leq d} v_i \partial_{u_i} H(t,u),
$$
$dS$ being the Lebesgue measure on $\mathbb{T}^{d-1}$. 

We say that that the solution $(\rho,{\boldsymbol p})$ has finite energy if its components belong to $L^2([0,T],H^1(D^d))$:
$$\int_0^T ds \left(\int_{D^d} \|\nabla \rho(s,u)\|^2 du\right)<\infty,$$
and 
$$\int_0^T ds \left(\int_{D^d} \|\nabla p_k(s,u)\|^2 du\right)<\infty,$$
for $k=1,\ldots,d$, where $\nabla f$ represents the generalized gradient of the function $f$.

In Section \ref{sec6} we prove that there exists at most one weak solution of the problem (\ref{problema2}).

\begin{theorem}\label{limhyd2} Let $(\mu^N)_N$ be a sequence of probability measures on $X_N$ associated to the profile $(\rho_0,{\boldsymbol p}_0)$. Then, for every $t\geq 0$, for every continuous function $H:D^d\to \mathbb{ R}$ vanishing at the boundary, and for every $\delta > 0$,
$$\lim_{N\to\infty} \mathbb{ P}_{\mu^N} \left[ \left| \frac{1}{N^d} \sum_{x\in D_N^d} H\left(\frac{x}{N}\right) I_0(\eta_x(t)) - \int_{D^d} H(u) {\rho}(t,u)du\right|>\delta\right] = 0,$$
and for $1\leq k \leq d$
$$\lim_{N\to\infty} \mathbb{ P}_{\mu^N} \left[ \left| \frac{1}{N^d} \sum_{x\in D_N^d} H\left(\frac{x}{N}\right) I_k(\eta_x(t)) - \int_{D^d} H(u) {p_k}(t,u)du\right|>\delta\right] = 0,$$
where $({\rho},{\boldsymbol p})$ has finite energy and is the unique weak solution of equation (\ref{problema2}).
\end{theorem}
The strategy to prove Theorem \ref{limhyd2} is to use a replacement lemma, together with some estimates
on Dirichlet forms and entropies for this boundary driven process.

\section{Hydrodynamic limit for the boundary driven process}\label{sec4}
Fix $T>0$, let $\mathcal{ M}_{+}$
be the space of finite positive measures on $D^d$ endowed with the weak
topology, and let $\mathcal{ M}$ be the space of bounded variation signed measures on $D^d$
endowed with the weak topology. Let $\mathcal{ M}_{+}\times\mathcal{ M}^{d}$ be the cartesian
product of these spaces endowed with the product topology, which
is metrizable.

Recall that the conserved quantities are the mass and momentum presented in subsection \ref{massmomentum}. For $k=0,\ldots,d$, denote by $\pi_t^{k,N}$ the empirical measure associated to the $k$th conserved quantity:
\begin{equation}\label{empmea}
\pi_t^{k,N} = \frac{1}{N^d} \sum_{x \in D_N^d} I_k(\eta_x(t)) \delta_{x/N},
\end{equation}
where $\delta_u$ stands for the Dirac measure supported on $u$. We denote by $<\pi_t^{k,N}, H>$ the integral of a test function $H$ with respect to an empirical measure $\pi_t^{k,N}$.

Let $D([0,T],\mathcal{ M}_{+}\times \mathcal{ M}^{d})$ be the set of right continuous functions with left limits taking values on $\mathcal{ M}_{+}\times\mathcal{ M}^{d}$. We consider the sequence of probability measures on $D([0,T],\mathcal{ M}_{+}\times \mathcal{ M}^{d})$ $(Q_N)_N$ that corresponds to the Markov process $\pi_t^{N} = (\pi_t^{0,N},\ldots,\pi_t^{d,N})$ starting from $\mu_N$.

Let $V$ be an open neighborhood of $D^d$, and consider, for each $v\in\mathcal{V}$, smooth functions $\kappa_k^v:V\to(0,1)$ in $C^2(V)$, for $k=0,\ldots,d$. We assume that each $\kappa_k^v$ has its image contained in some compact subset of $(0,1)$, that the restriction of $\kappa = \sum_{v\in\mathcal{V}}(\kappa_0^v,v_1\kappa_1^v,\ldots,v_d\kappa_d^v)$ to $\{0\}\times\mathbb{T}^{d-1}$ equals the vector valued function $a(\cdot)$ defined in \eqref{funcaod}, and that the restriction of $\kappa$ to $\{1\}\times\mathbb{T}^{d-1}$ equals the vector valued function $b(\cdot)$, also defined in \eqref{funcaod}, in the sense that $\kappa(x) = d(x_1,\tilde{x})$ if $x\in\{0,1\}\times\mathbb{T}^{d-1}$.

Further, we
may choose $\kappa$ for which there exists a constant $\theta>0$ such
that:
\begin{eqnarray*}
\kappa(u_1,\tilde{u}) =  d(-1,\tilde{u}) & 
\qquad\hbox{ if } \; 0\leq u_1\leq \theta\, , \\ 
\kappa(u_1,\tilde{u}) =  d(1,\tilde{u})\;\;\, & 
\,\hbox{ if } \;1-\theta\leq u_1\leq 1\, ,
\end{eqnarray*}
for all $\tilde{u}\in\bb T^{d-1}$. In that case, for every $N$ large
enough, $\nu_{\kappa}^N$ is reversible for the process with
generator $\mc L_{N}^b$ and then $\langle -N^2\mc L_{N}^b f,
f\rangle_{\nu^N_{\kappa}}$ is positive.

We then consider $\nu_\kappa^N$ the product measure on $X_N$ with marginals given by 
$$\nu_\kappa^N\{\eta: \eta(x,\cdot) = \xi\} = m_{\Lambda(\kappa(x))}(\xi),$$
where $m_\lambda(\cdot)$ was defined in \eqref{mprod}. Note that with this choice, for $N$ sufficiently large, we have that if $x\in\{1\}\times \mathbb{T}_N^{d-1}$, then $E_{\nu_\kappa^N}[\eta(x,v)] = \alpha_v(\tilde{x}/N)$ and if $x\in\{N-1\}\times\mathbb{T}_N^{d-1}$, then $E_{\nu_\kappa^N}[\eta(x,v)] = \beta_v(\tilde{x}/N).$

\subsection{Entropy estimates}
Let us recall some definitions. Recall that ${S}_t^N$ is the
semigroup associated to the generator ${\mathcal{L}}_N = N^2(\mathcal{L}_N^{ex} + \mathcal{L}_N^{c} +\mathcal{L}_N^{b})$. Denote by $f_t = f_t^N$
the Radon-Nikodym derivative of $\mu^N {S}_t^N$ with respect to $\nu_{\kappa}^N$. For each function $f:X_N\to\bb R$, let $D_{\nu_{\kappa}^N}(f)$ be
$$D_{\nu_{\kappa}^N}(f) = D_{\nu_{\kappa}^N}^{ex}(f) + D_{\nu_{\kappa}^N}^c(f) + D_{\nu_{\kappa}^N}^b(f),$$
where
$$D_{\nu_{\kappa}^N}^{ex}(f) = \sum_{v\in\mathcal{V}}\sum_{x\in D_N^d} \sum_{x+z\in D_N^d}P_N(z-x,v)\int \left[\sqrt{f(\eta^{x,z,v})} - \sqrt{f(\eta)}\right]^2 \nu_\kappa^n(d\eta),$$
$$D_{\nu_{\kappa}^N}^{c}(f) = \sum_{q\in\mathcal{Q}}\sum_{x\in D_N^d} \int p(x,q,\eta)\left[\sqrt{f(\eta^{x,q})} - \sqrt{f(\eta)} \right]^2\nu_\kappa^N(d\eta),$$
and
\begin{align*}
D_{\nu_{\kappa}^N}^{b}(f) = \sum_{v\in\mathcal{V}}\sum_{x\in\{1\}\times\mathbb{T}_N^{d-1}} \int &[\alpha_v(\tilde{x}/N)(1-\eta(x,v)) + (1-\alpha_v(\tilde{x}/N))\eta(x,v)]\times\\
&\times\left[\sqrt{f(\sigma^{x,v}\eta)}-\sqrt{f(\eta)}\right]^2\nu_\kappa^N(d\eta) \; +\\
+\;\sum_{v\in\mathcal{V}}\sum_{x\in\{N-1\}\times\mathbb{T}_N^{d-1}}\int &[\beta_v(\tilde{x}/N)(1-\eta(x,v)) + (1-\beta_v(\tilde{x}/N))\eta(x,v)]\times\\
&\times\left[\sqrt{f(\sigma^{x,v}\eta)}-\sqrt{f(\eta)}\right]^2\nu_\kappa^N(d\eta).
\end{align*}
\begin{proposition} There exists a finite constant $C = C(\alpha,\beta)$ such that
\begin{equation}\label{entropy}
\partial_t H(\mu^N {S}_t^N\vert \nu_{\kappa}^N) \leq -N^2 D_{\nu_{\kappa}^N}(f_t) + CN^d.
\end{equation}
\end{proposition}
\emph{Proof}: Denote by $\mathcal{L}_\nu^\ast$ the adjoint operator of ${\mathcal{L}}_N$ with respect to $\nu_{\kappa}^N$. Then,
$f_t$ is the solution of the forward equation
$$\left\{
\begin{array}{c}
\partial_t f_t = N^2 \mathcal{L}_\nu^\ast f_t,\\
f_0 = d\mu^N/d\nu_{\kappa}^N.
\end{array}\right.
$$
Thus,
\begin{eqnarray*}
\partial_t H(\mu^N {S}_t^N\vert \nu_{\kappa}^N)\!\!\!\! &=&\!\!\!\!\!\! \int N^2 \mathcal{L}_\nu^\ast f_t \log f_t d\nu_{\kappa}^N+\!\!\!\int N^2 {\mathcal{L}}_\nu^\ast f_t d\nu_{\kappa}^N =\!\!\! \int f_t N^2\mathcal{L}_N\log f_t d\nu_{\kappa}^N\\
&=& N^2 \int f_t({\mathcal{L}}_N \log f_t - \frac{{\mathcal{L}}_N f_t}{f_t})d\nu_{\kappa}^N + N^2 \int {\mathcal{L}}_N f_t d\nu_{\kappa}^N.
\end{eqnarray*}
Note that the last term is the price paid for not using an invariant measure.

Since for every $a,b>0$, $a\log(b/a) - (b-a)$ is less than or equal to $-(\sqrt{b}-\sqrt{a})^2$, for every $x,y\in D_N^d$, we have
$$f_t \mathcal{L}_{x,y,v}^{ex} \log f_t - \mathcal{L}_{x,y,v}^{ex} f_t \leq - P_N(y-x,v)\left[\sqrt{f_t(\eta^{x,y,v})}-\sqrt{f_t(\eta)}\right]^2.$$

An analogous calculation for the other parts of the generator permits to conclude that
$$N^2 \int f_t({\mathcal{L}}_N \log f_t - \frac{{\mathcal{L}}_N f_t}{f_t})d\nu_{\kappa}^N \leq -N^2 D_{\nu_{\kappa}^N}(f_t).$$
To conclude the proposition we need a bound for $N^2\int {\mathcal{L}}_N f_t d\nu_{\kappa}^N$. Let us write it explicitly:
$$N^2\int{\mathcal{L}}_N f_t d\nu_{\kappa}^N = N^2 \int (\mathcal{L}_N^{ex,1} f_t + \mathcal{L}_N^{ex,2} f_t + \mathcal{L}_N^c f_t + \mathcal{L}_N^b f_t)d\nu_{\kappa}^N.$$
Now, we compute each term inside this integral separately.
\begin{eqnarray*}
N^2 \int \mathcal{L}_N^{ex,1} f_t d\nu_{\kappa}^N \!\!\!&=&\!\!\!
N^2 \int \sum_{v\in\mathcal{V}}\sum_{x\in D_N^d}\sum_{j=1}^d [f(\eta - \mf{d}_{x,v}+\mf{d}_{x+e_j,v})-f(\eta)]d\nu_{\kappa}^N\\
\!\!\!&+&\!\!\! N^2 \int \sum_{v\in\mathcal{V}}\sum_{x\in D_N^d}\sum_{j=1}^d [f(\eta - \mf{d}_{x,v}+\mf{d}_{x-e_j,v})-f(\eta)]d\nu_{\kappa}^N,
\end{eqnarray*}
where $\mf{d}_{x,v}$ represents a configuration with one particle at position $x$ and velocity $v$, and no particles elsewhere. Then, if we let
$$\gamma_{x,v} = \theta_v(\Lambda(\kappa(x)))/(1-\theta_v(\Lambda(\kappa(x)))),$$
the change of variables
$\eta - \mf{d}_{x,v} + \mf{d}_{x+e_j,v} = \xi$, changes the measure as $d\nu_{\kappa}^N(\eta)/d\nu_{\kappa}^N(\xi) = \gamma_{x,v}/\gamma_{x+e_j,v}$. Hence, after changing
the variables, we obtain
\begin{eqnarray*}
N^2 \int \mathcal{L}_N^{ex,1} f_t d\nu_{\kappa}^N &=& N^2 \sum_{v\in\mathcal{V}}\sum_{j=1}^d\int\sum_{x\in D_N^d}\left[\frac{\gamma_{x,v}}{\gamma_{x+e_j,v}}-1\right]f_t(\eta)d\nu_{\kappa}^N\\
&+& N^2 \sum_{v\in\mathcal{V}}\sum_{j=1}^d\int\sum_{x\in D_N^d}\left[\frac{\gamma_{x,v}}{\gamma_{x-e_j,v}}-1\right]f_t(\eta)d\nu_{\kappa}^N\\
&=& \sum_{v\in\mathcal{V}}\sum_{j=1}^d\int\sum_{x\in D_N^d}\frac{\Delta_N \gamma (x,v)}{\gamma_{x,v}} f_t(\eta) d\nu_{\kappa}^N\\
&+& N\sum_{v\in\mathcal{V}}\int\sum_{\substack{x\in D_N^d\\x_1=1}}\frac{\partial_{u_1}^N\gamma(x,v)}{\gamma_{x,v}}f_t(\eta)d\nu_{\kappa}^N\\
&-& N\sum_{v\in\mathcal{V}}\int\sum_{\substack{x\in D_N^d\\x_1=N-1}}\frac{\partial_{u_1}^N\gamma(x,v)}{\gamma_{x,v}}f_t(\eta)d\nu_{\kappa}^N.
\end{eqnarray*}

Since $\gamma_{x,v}$ is smooth and does not vanish, we can bound the above quantity by $C_1 N^d$, where $C_1$ is a constant depending only on $\alpha$ and $\beta$. By a similar approach, one may conclude that
$$N^2 \int \mathcal{L}_N^{ex,2} f_t d\nu_{\kappa}^N \leq \sum_{v\in\mathcal{V}}\sum_{j=1}^d v_j \sum_{x\in D_N^d} \frac{\partial_{u_i}^N \gamma(x,v)}{\gamma_{x,v}},$$
which is clearly bounded by $C_2 N^d$, where $C_2$ is a constant depending only on $\alpha$ and $\beta$.

We now move to the generator with respect to collision. The change of variables $\eta^{y,q} = \xi$ changes the measure as $d\nu_{\kappa}^N(\eta)/d\nu_{\kappa}^N(\xi) = (\gamma_{y,v}\gamma_{y,w})/(\gamma_{y,v'}\gamma_{y,w'})$,
where $v+w = v'+w'$. Then, clearly, $(\gamma_{y,v}\gamma_{y,w})/(\gamma_{y,v'}\gamma_{y,w'}) = 1$, and therefore
$$N^2 \int \mathcal{L}_N^{c} f_t d\nu_{\kappa}^N = 0.$$

Lastly, we note that the change of variables $\sigma^{x,v} \eta = \xi$ changes the measure $d\nu_{\kappa}^N(\eta)/d\nu_{\kappa}^N(\xi) = \alpha_v(\tilde{x}/N)/(1-\alpha_v(\tilde{x}/N))$ or $(1-\alpha_v(\tilde{x}/N))/\alpha_v(\tilde{x}/N)$,
depending on whether there is or there is not a particle at the site $x$ with velocity $v$, and analogously for $\beta$. Therefore, a simple computation shows that
$$N^2 \int \mathcal{L}_N^{b} f_t d\nu_{\kappa}^N = 0.$$
which concludes the Proposition. $\square$\\
Let $<f,g>_\nu$ be the inner product in $L^2(\nu)$ of $f$ and $g$: 
$$<f,g>_\nu = \int fg d\nu.$$
\begin{proposition}\label{dirichlet}
There exist constants $C_1>0$ and $C_2 = C_2(\alpha,\beta)>0$ such that for every density $f$ with respect to $\nu_\kappa^N$, then
$$<{\mathcal{ L}}_N\sqrt{f},\sqrt{f}>_{\nu_\kappa^N} \leq -C_1 D_{\nu_\kappa^N}(f) + C_2 N^{d-2}.$$
\end{proposition}
\emph{Proof:} A simple computation permits to conclude that $D_{\nu_\kappa^N}^c$ and $D_{\nu_\kappa^N}^b$ are both non-negative. Finally, the computation for $D_{\nu_\kappa^N}^{ex}$ follows the same lines as those on the proof of Proposition \ref{entropy}, and on Lemmas \ref{lemabound1} and \ref{lemabound2}, and is therefore omitted. $\square$

\subsection{Replacement lemma for the boundary}
Fix $k=0,\ldots,d$, a continuous function $G:[0,T]\times\mathbb{T}^{d-1}\to\mathbb{R}^{d+1}$, and consider the quantities
\begin{equation*}
V_k^{-}(s,\eta,\alpha,G) = \frac{1}{N^{d-1}}\sum_{\tilde{x}\in\mathbb{T}_N^{d-1}} G_k(s,\tilde{x}/N) \Big(I_k(\eta_{(1,\tilde{x})}(s))- \sum_{v\in\mathcal{V}}v_k\alpha_v(\tilde{x}/N)\Big),
\end{equation*}
\begin{equation*}
V_k^{+}(s,\eta,\beta,G) = \frac{1}{N^{d-1}}\sum_{\tilde{x}\in\mathbb{T}_N^{d-1}} G_k(s,\tilde{x}/N) \Big(I_k(\eta_{(N-1,\tilde{x})}(s))- \sum_{v\in\mathcal{V}}v_k\beta_v(\tilde{x}/N)\Big),
\end{equation*}
\begin{equation*}
V_k^2(s,\eta,\alpha,G) = \frac{1}{N^{d-1}}\sum_{\tilde{x}\in\mathbb{T}_N^{d-1}} G_k(s,\tilde{x})\Big(I_k(\eta_{(1,\tilde{x})}(s))
- \frac{1}{N\epsilon}\sum_{x_1=1}^{N\epsilon-1} I_k(\eta_{(1,\tilde{x})}(s))\Big),
\end{equation*}
and
\begin{align*}
V_k^2(s,\eta,\beta,G) =\frac{1}{N^{d-1}}\sum_{\tilde{x}\in\mathbb{T}_N^{d-1}} G_k(s,\tilde{x})\Big(I_k(\eta_{(N-1,\tilde{x})}(s))
- \frac{1}{N\epsilon}\sum_{x_1=N(1-\epsilon)-1}^{N-1} I_k(\eta_{(N-1,\tilde{x})}(s))\Big),
\end{align*}
where $s\in [0,T]$, and $G_k$, $0 \leq k\leq d$ are the components of function $G$.\\

The main result of this subsection is the following Lemma:

\begin{lemma}\label{lemabound3}
For each $0\leq t\leq T$, $0\leq k \leq d$, and $G:[0,T]\times D^d\to\mathbb{R}$ continuous,
\begin{align*}
\limsup_{N\to\infty} E_{\mu^N}\left[\left\vert\int_0^t ds V_k^j(s,\eta,\zeta,G)\right\vert\right] = 0,
\end{align*}
where $j=1,2$, and $\zeta=\alpha,\beta$.
\end{lemma}
\emph{Proof:}  It is clear that $V_k^j$ is bounded for each $0\leq k\leq d$, and $j=1,2$. By the entropy inequality,
\begin{align*}
E_{\mu^N} \Big[ \Big\vert \int_0^t ds &V_k^j(s,\eta,\zeta,G)\Big\vert \Big]\leq\\
\leq \frac{H(\mu^N|\nu_{\kappa}^N)}{AN^d}& + \frac{1}{AN^d} \log E_{\nu_{\kappa}^N} \left[\exp\left\{ \left\vert \int_0^t ds AN^d V_k^j(s,\eta,\zeta,G)\right\vert\right\}\right],
\end{align*}
for all $A>0$. We have that the first term on the right-hand side is bounded by $CA^{-1}$, for some constant $C$. To prove this result we must show that
the limit of the second term is less than or equal to 0 as $N\to\infty$ for some suitable choice of $A>0$. Since $e^{|x|}\leq e^x + e^{-x}$ and $\limsup_{N\to\infty} N^{-d}\log\{a_N + b_N\}$ $\leq$ $\max\{\limsup_{N\to\infty}N^{-d}\log(a_N),$ $\limsup_{N\to\infty}N^{-d}\log(b_N)\}$,
replacing $V_k^j$ by $-V_k^j$, or more precisely, replacing $G_k$ by $-G_k$, we are able to conclude that we only need to prove the previous statement without the
absolute values in the exponent. Let $W_k(s) = AN^d V_k^j(s,\eta,\zeta,G)$. Then, by Feynman-Kac's formula (see, for instance, \cite{bkl,KL}), we have
$$E_{\nu_\kappa^N}\left[\exp\left\{ \int_0^t ds AN^d V_k^j(s,\eta,\zeta,G)\right\} \right] = <S_{0,t}^{W_k} 1,1>_{\nu_\kappa^N},$$
where $S_{s,t}^{W_k}$ is a semigroup associated to the operator ${\mathcal{L}}_t^W = {\mathcal{L}} + W_k(t)$, for more details see \cite[A.1.7]{KL}, see also \cite{bkl}. Then, by Cauchy-Schwarz
$$<S_{0,t}^{W_k}1,1>_{\nu_\kappa^N} \leq <S_{0,t}^{W_k}1,S_{0,t}^{W_k}1>_{\nu_\kappa^N}^{1/2}.$$
On the other hand, since $W_k$ is bounded, the adjoint in $L^2(\nu_{\kappa}^N)$ of ${\mathcal{L}}_t^W$, $L_t^{W,\ast}$, is equal to ${\mathcal{L}}_N^\ast+W_k(t)$. We have that
\begin{eqnarray*}
\partial_s <S_{s,t}^{W_k}1, S_{s,t}^{W_k}1>_{\nu_{\kappa}^N}&=& <({\mathcal{L}}_t^{W_k}+{\mathcal{L}}_t^{W_k,\ast})S_{s,t}^{W_k}1,S_{s,t}^{W_k} 1>\\
&=& 2<{\mathcal{L}}_t^{W_k}S_{s,t}^{W_k}1,S_{s,t}^{W_k} 1>\leq \lambda_{W_k}(s) <S_{s,t}^{W_k}1,S_{s,t}^{W_k}1>_{\nu_\kappa^N},
\end{eqnarray*}
where $\lambda_{W_k}(s) = \sup_{\|f\|_{L^2(\nu_\kappa^N)=1}} \left\{ <W_k(s), f>_{\nu_\kappa^N} + <{\mathcal{L}}_N f,f>_{\nu_\kappa^N}\right\}$.
Therefore, we obtained that
\begin{align*}
\frac{1}{AN^d} \log &E_{\nu_{\kappa}^N} \Big[\exp\Big\{ \Big\vert \int_0^t ds AN^d V_k^j(s,\eta,\zeta,G)\Big\vert\Big\}\Big]\leq\\
&\leq \int_0^t ds \sup_f \left\{\int V_k^j(s,\eta,\zeta,G) f(\eta(s)) d\nu_{\kappa}^N  + \frac{<{\mathcal{L}}_N\sqrt{f},\sqrt{f}>_{\nu_\kappa^N}}{AN^{d-2}}\right\}.
\end{align*}
In this formula the supremum is taken over all densities $f$ with respect to $\nu_{\kappa}^N$, and recall that $<f,g>_\nu$ stands for the inner product in $L^2(\nu)$ of $f$ and $g$. An application of Proposition \ref{dirichlet} permits to conclude that $<{\mathcal{L}}_N\sqrt{f},\sqrt{f}>_{\nu_\kappa^N}$ is bounded above by $C N^{d-2}$, where $C>0$ is some constant. Thus, if we choose, for instance, $A=N$, the proof follows from an application of the auxiliary Lemmas \ref{lemabound1} and \ref{lemabound2} given below. $\square$

\begin{lemma}\label{lemabound1} 
For every $0\leq t\leq T,$ $0\leq k \leq d$, and every continuous $G:[0,T]\times\mathbb{T}^{d-1}\to\mathbb{R}^{d+1}$,
\begin{align*}
\limsup_{N\to\infty} E_{\mu^N}\left[\int_0^t dsV_k^1(s,\eta,\zeta,G)\right] = 0,
\end{align*}
where $\zeta=\alpha,\beta$.
\end{lemma}
\emph{Proof:} We will only prove for $\alpha$, since for $\beta$ the proof is entirely analogous.
Note that $G$ is continuous and its domain is compact, hence, we may prove the above result without $G$. Set $\overline{f}_t = 1/t \int_0^t f_s ds$. With this notation
we can write the expectation above, without $G$, as
\begin{eqnarray*}
\frac{t}{N^{d-1}}\sum_{\tilde{x}\in\mathbb{T}_N^{d-1}}\int \overline{f}_t(\eta) \left[ I_k(\eta_{(1,\tilde{x})} - \sum_{v\in\mathcal{V}}v_k\alpha_v(\tilde{x}/N)\right]d\nu_{\kappa}^N\\
= \frac{t}{N^{d-1}}\sum_{\tilde{x}\in\mathbb{T}_N^{d-1}}\sum_{v\in\mathcal{V}}v_k \int \overline{f}_t(\eta) \left[\eta((1,\tilde{x}),v)-\alpha_v(\tilde{x}/N)\right] d\nu_{\kappa}^N.
\end{eqnarray*}
Then, splitting the integral into the integral over the sets $[\eta((1,\tilde{x}),v) = 0]$ and $[\eta((1,\tilde{x}),v)=1]$, and changing the variables as $1-\eta(x_N,v) = \xi$,
we obtain
\begin{eqnarray*}
\frac{t}{N^{d-1}}\sum_{\tilde{x}\in\mathbb{T}_N^{d-1}}\int \overline{f}_t(\eta) \left[ I_k(\eta_{(1,\tilde{x})} - \sum_{v\in\mathcal{V}}v_k\alpha_v(\tilde{x}/N)\right]d\nu_{\kappa}^N\\
= \frac{t}{N^{d-1}}\sum_{\tilde{x}\in\mathbb{T}_N^{d-1}}\sum_{v\in\mathcal{V}}v_k \int P_{\alpha,\eta}\left[\overline{f}_t(\eta) - \overline{f}_t(\eta-\mf{d}_{(1,\tilde{x}),v})\right] d\nu_{\kappa}^N,
\end{eqnarray*}
where
$$P_{\alpha,\eta} = \alpha_v(\tilde{x}/N)(1-\eta((1,\tilde{x}),v))+(1-\alpha_v(\tilde{x}/N))\eta((1,\tilde{x}),v).$$
Writing $\{a-b\} = \{\overline{f}_t(\eta) - \overline{f}_t(\eta-\mf{d}_{(1,\tilde{x}),v})\}$ as $\{\sqrt{a}-\sqrt{b}\}\{\sqrt{a}+\sqrt{b}\}$ and applying
Cauchy-Schwarz, the above expression is bounded by
$$\frac{2 t \sum_{v\in\mathcal{V}} v_k}{A} + \frac{t}{N^{d-1}}AD_{\nu_{\kappa}^N,b}(\overline{f}_t),$$
where $D_{\nu_{\kappa}^N,b}(\overline{f}_t)$ is the Dirichlet form of $\overline{f}_t$ with respect to $\mathcal{L}_N^b$.
Then, choosing $A=\sqrt{N}$, the proof of the Lemma follows from an application of Proposition \ref{entropy} together with the fact that the Dirichlet form
is convex. $\square$\\\\
The next Lemma concludes the boundary behavior of the particle system.
\begin{lemma}\label{lemabound2}
For each $0\leq t\leq T$, $0\leq k \leq d$, and continuous $G:[0,T]\times D^d$,
\begin{align*}
\limsup_{\epsilon\to 0}\limsup_{N\to\infty} E_\mu^N\left[ \int_0^t ds V_k^2(s,\eta,\zeta,G)\right]=0,
\end{align*}
where $\zeta=\alpha,\beta$.
\end{lemma}
\emph{Proof:} First of all, note that since $G$ is continuous and its domain $[0,T]\times D^d$ is compact, it is enough to prove
the result without the multiplying factor $G$. Moreover, we will only prove the first limit above, since the proof of the second one is entirely analogous.
Considering the notation used to prove Lemma \ref{lemabound1}, we may write the expectation above, without $G$, as
$$\frac{t}{N^{d-1}}\sum_{\tilde{x}\in\mathbb{T}_N^{d-1}} \int \left[ I_k (\eta_{(1,\tilde{x})}) - \frac{1}{N\epsilon}\sum_{x_1=1}^{N\epsilon-1} I_k(\eta_{(x_1,\tilde{x})})\right] d\nu_{\kappa}^N.$$
We now obtain, by a change of variables and a telescopic sum, that the absolute value of the above expression is bounded above by
$$\left\vert\frac{t}{N^{d-1}}\sum_{\tilde{x}\in\mathbb{T}_N^{d-1}} \frac{1}{N\epsilon}\sum_{y=1}^{N\epsilon-1} \sum_{x_1=1}^{y-1}K_1\int\left[\overline{f}_t(\prod_{i=1}^{x_1}\tau_{z_i}(\eta))-\overline{f}_t(\prod_{i=1}^{x_1-1}\tau_{z_i}(\eta))\right] d\nu_{\kappa}^N.\right\vert,$$
where $K_1$ is a constant which depends on $\alpha$, $\beta$ and $d$, $z_1=1,\ldots,z_{y-1}=y$ is the path from the origin to $y$ across the first coordinate of the space,
and $\tau_{z_1}(\eta)\cdots\tau_{z_i}(\eta)$ is the sequence of nearest neighbor exchanges that represents the path along $z_1,\ldots,z_i$. By Cauchy-Schwarz,
this expression is bounded above by
$$\frac{tA}{N^{d-1}} \sum_{\tilde{x}\in\mathbb{T}_N^{d-1}} \frac{1}{N\epsilon} \sum_{y=1}^{N\epsilon-1}\sum_{x_1=1}^{y-1} K_1 \int \left[ \sqrt{\overline{f}_t(\prod_{i=1}^{x_1}\tau_{z_i}(\eta))} - \sqrt{\overline{f}_t(\prod_{i=1}^{x_1-1}\tau_{z_i}(\eta))}\right]^2 d\nu_{\kappa}^N\; +$$
$$+\;\frac{t}{A N^{d-1}} \sum_{\tilde{x}\in\mathbb{T}_N^{d-1}} \frac{1}{N\epsilon} \sum_{y=1}^{N\epsilon-1}\sum_{x_1=1}^{y-1} K_1 \int \left[ \overline{f}_t(\prod_{i=1}^{x_1}\tau_{z_i}(\eta)) - \overline{f}_t(\prod_{i=1}^{x_1-1}\tau_{z_i}(\eta))\right] d\nu_{\kappa}^N,$$
for every $A>0$. Now, we can bound above the last expression by
$$\frac{tA K_1}{N^{d-1}} D_{\nu_{\kappa}^N}^{ex} (\overline{f}_t) + \frac{tK_2 N\epsilon}{A},$$
for every $A>0$, where $K_2$ is a constant that depends on $K_1$. Then, choosing $A=\sqrt{\epsilon}N$ and applying Proposition \ref{entropy}, we conclude
the proof of this Lemma. $\square$\\\\

\subsection{Tightness}
To prove tightness of the sequence $(Q_N)_N$, it is enough to prove that for every $k=0,\ldots,d$
\begin{align*}
\lim_{\delta\to 0} \limsup_{N\to\infty} \mathbb{ E}_{\mu^N}\left[ \sup_{|t-s|<\delta} \left| \frac{1}{N^d}\sum_{x\in D_N^d} H\left(\frac{x}{N}\right)\right.\right.I_k(\eta_x(t))
- \left.\left.\frac{1}{N^d} \sum_{x\in D_N^d} H\left(\frac{x}{N}\right)I_k(\eta_x(s))\right|\right] = 0,
\end{align*}
for any smooth test function $H:D^d\to\mathbb{R}$ vanishing at the boundary. 

Fix $0\leq k \leq d$, then, by Dynkin's formula
\begin{equation}\label{dynkin}
M_t^k = <\pi_t^{k,N}, H> - <\pi_0^{k,N}> - \int_0^t \mathcal{ L}_N <\pi_s^{k,N},H> ds
\end{equation}
is a martingale. On the other hand,
$$\mathbb{ E}_{\mu^N} [M_t^k]^2 = \mathbb{ E}_{\mu^N} \left[ \int_0^t \left\{\mathcal{ L}_N <\pi_s^{k,N},H>^2 - 2<\pi_s^{k,N},H>\mathcal{ L}_N <\pi_s^{k,N},H> \right\}ds\right].$$
Writing the above expression as four sums, the first corresponds to the nearest neighbor symmetric exclusion process and the other corresponds to the asymmetric exclusion process, the third and fourth corresponding to the collision and boundary parts of the dynamics, respectively. A long, albeit simple computation shows that all of these sums are of order $\mathcal{ O}(N^{-d})$, and therefore, the right-hand side of the above expression is of the same order. Thus, by Doob's inequality, $\mathbb{ E}_{\mu^N} [\sup_{0\leq s\leq t} (M_t^k)^2] = \mathcal{ O}(N^{-d})$. 

Hence, by \eqref{dynkin} and the above estimates, we have
$$\frac{1}{N^d} \sum_{x\in D_N^d} H\left(\frac{x}{N}\right) I_k(\eta_x(t)) = \frac{1}{N^d} \sum_{x\in D_N^d} H\left(\frac{x}{N}\right) I_k(\eta_x(s))\;+$$
$$+ \frac{1}{N^d} \sum_{j=1}^d\sum_{x,z\in D_N^d}\sum_{v\in\mathcal{ V}} \int_s^t p(z,v) v_k \eta_r(0,v)[1-\eta_r(z,v)]z_j (\partial_{u_j} H)\left(\frac{x}{N}\right) dr\;+$$
$$+\frac{1}{2N^d} \sum_{x\in D_N^d} \int_s^t (\Delta H)\left(\frac{x}{N}\right) I_k(\eta_x(r))dr+\frac{1}{N^{d-1}}\sum_{\substack{x \in D_N^d\\x_1=N-1}}  \int_s^t\partial_{u_1} H\left(\frac{x}{N}\right)I_k(\eta_x(r))dr$$
$$-\frac{1}{N^{d-1}}\sum_{\substack{x \in D_N^d\\x_1=1}}  \int_s^t\partial_{u_1} H\left(\frac{x}{N}\right)I_k(\eta_x(r))dr + R_N + \mathcal{ O}(N^{-d})+\mathcal{ O}(N^{-1}),$$
where the terms were obtained from $\mathcal{ L}_N <\pi_s^{k,N},H>$ by means of summation by parts, and the replacement of discrete derivatives and discrete Laplacian by the continuous ones, and $R_N$ is the error coming from such replacements. Since $p$ is of finite range, the error $R_N$ is uniformly of order $\mathcal{ O}(N^{-1})$. Finally, by using Lemma \ref{lemabound3} and a calculation similar to the one found in equation \eqref{contaboundary}, we have that $\mathcal{ L}_N^b <\pi_s^{k,N},H>=\mathcal{O}(N^{-1})$. Tightness thus follows from the above estimates.

Our next goal is to prove the replacement lemma. To do so, we need the following result known as equivalence of ensembles, which will be used in the proofs of the one block estimate and of the two block estimate.
\subsection{Equivalence of ensembles}\label{ensemble}
Fix $L\geq 1$ and a configuration $\eta$, let ${\boldsymbol I}^L(x,\eta):={\boldsymbol I}^L(x) = (I_0^L(x),\ldots,I_d^L(x))$ be the average of the conserved quantities in a cube of the length $L$ centered at $x$:
$${\boldsymbol I}^L(x)= \frac{1}{|\Lambda_L|} \sum_{z\in x+\Lambda_L} {\boldsymbol I}(\eta_z),$$
where, $\Lambda_L = \{-L,\ldots,L\}^d$ and $|\Lambda_L| = (2L+1)^d$ is the discrete volume of box $\Lambda_L$.

Let ${\mf V}_L$ be the set of all possible values of ${\boldsymbol I}^L(0,\eta)$ when $\eta$ runs over $\left(\{0,1\}^\mathcal{ V}\right)^{\Lambda_L}$, that is, 
$${\mf V}_L = \left\{{\boldsymbol I}^L(0,\eta); \eta \in \left(\{0,1\}^\mathcal{ V}\right)^{\Lambda_L}\right\}.$$

Note that ${\mf V}_L$ is a finite subset of the convex envelope of $\left\{ {\boldsymbol I}(\xi): \xi \in \{0,1\}^\mathcal{ V}\right\}$. The set of configurations $\left( \{0,1\}^\mathcal{ V}\right)^{\Lambda_L}$ splits into invariant subsets: for ${\boldsymbol i}$ in ${\mf V}_L$, let
$$\mathcal{ H}_L({\boldsymbol i}) := \left\{\eta \in \left( \{0,1\}^\mathcal{ V}\right)^{\Lambda_L}: {\boldsymbol I}^L(0) = {\boldsymbol i}\right\}.$$
For each ${\boldsymbol i}$ in ${\mf V}_L$, define the canonical measure $\nu_{L,{\boldsymbol i}}$ as the uniform probability measure on $\mathcal{ H}_L({\boldsymbol i})$. Note that for every $\bs \lambda$ in $\mathbb{ R}^{d+1}$\label{measunif}
$$\nu_{\Lambda_L,{\boldsymbol i}} (\cdot ) = \mu^{\Lambda_L}_{{\boldsymbol \lambda}} \left(  \cdot \left\vert {\boldsymbol I}^L ={\boldsymbol i}\right.\right).$$ 
Let $< g; f>_{\mu}$ stands for the covariance of $g$ and $f$ with respect to $\mu$: $< g; f>_{\mu} = E_{\mu}[fg] - E_{\mu}[f] E_{\mu}[g]$.
\begin{proposition} \label{ensembles}(Equivalence of ensembles): Fix a cube $\Lambda_\ell\subset\Lambda_L$. For each ${\boldsymbol i}\in {\mf V}_L$, denote by $\nu^{\ell}$ the projection of the canonical measure $\nu_{\Lambda_L,{\boldsymbol i}}$ on $\Lambda_\ell$ and by $\mu^{\ell}$ the projection of the grand canonical measure $\mu^L_{{\boldsymbol \Lambda} ({\boldsymbol i})}$ on $\Lambda_\ell$. Then, there exists a finite constant $C(\ell, \mathcal V)$, depending only on $\ell$ and $\mathcal V$, such that 
$$\left|E_{\mu^{\ell}}[f]-E_{\nu^{\ell}}[f]\right| \leq \frac{C(\ell, \mathcal{ V})}{|\Lambda_L|} < f; f>^{1/2}_{\mu^{\ell}}$$
for every function $f:\left(\{0,1\}^\mathcal{ V}\right)^{\Lambda_{\ell}}\mapsto \mathbb{ R}$.
\end{proposition}
The proof of Proposition \ref{ensembles} can be found in Beltr\'an and Landim \cite{BL}.
\subsection{Replacement lemma}\label{subreplemma}
We now state the replacement lemma that will allow us to prove that the limit points $Q$ are concentrated on weak solutions of (\ref{problema2}).
\begin{lemma}\label{replacement}
(Replacement lemma): For all $\delta>0$, $1\leq j\leq d$, $0\leq k\leq d$: 
$$\limsup_{\epsilon\to 0}\limsup_{N\to\infty} \mathbb{ P}_{\mu^N} \left[ \int_0^T \frac{1}{N^d} \sum_{x\in D_N^d} \tau_x V_{\epsilon N}^{j,k} (\eta(s))ds \geq \delta \right] = 0,$$
where 
\begin{align}
V_\ell^{j,k}(\eta) = \left| \frac{1}{(2\ell +1)^d}\sum_{y\in{\Lambda_\ell}} \sum_{v\in\mathcal{ V}} v_k \sum_{z\in\mathbb{ Z}^d} p(z,v)z_j\right.\tau_y( \eta(0,v)[1-\eta(z,v)] ) 
- \left.\sum_{v\in\mathcal{ V}} v_j v_k \chi(\theta_v(\Lambda({\boldsymbol I}^\ell(0))))\right|.\label{Vrepl}
\end{align}
\end{lemma}

Note that $V_{\epsilon N}^{j,k}$ is well-defined for large $N$ since $p(\cdot,v)$ is of finite range. 
We now observe that Propositions \ref{entropy} and \ref{dirichlet} permit us to prove the following replacement lemma for the boundary driven exclusion process by using the process without the boundary part of the generator (see \cite{lms} for further details). We postpone the rest of the proof to Section \ref{sec5}. 

\subsection{Energy estimates}
We will now define some quantities to prove that each component of the solution vector belongs, in fact, to $H^1([0,T]\times D^d)$. The proof is similar to the one found in \cite{KLO}.

Let the energy $\mathcal{ Q}:D([0,T],\mathcal{ M})\to[0,\infty]$ be given by
$$\mathcal{ Q}(\pi) = \sum_{i=1}^d \mathcal{ Q}_i(\pi),$$
with 
$$\mathcal{ Q}_i(\pi) = \sup_{G\in C_c^\infty(\Omega_T)}\left\{2\int_0^T dt <\pi_t,\partial_{u_i} G_t> - \int_0^T\int_{D^d} du G(t,u)^2\right\},$$
where $\Omega_T = (0,T)\times  D^d$ and $C_c^\infty(\Omega_T)$ stands for the set of infinitely differentiable functions (with respect to both the time and space) with compact support in $\Omega_T$. Let now, for any $G\in C_c^\infty(\Omega_T)$, $1\leq i\leq d$ and $C\geq 0$, 
$\mathcal{ Q}_{i,C}^G:D([0,T],\mathcal{ M})\to \mathbb{ R}$ be the functional given by
$$\mathcal{ Q}_{i,C}^G(\pi) = \int_0^T ds <\pi_s,\partial_{u_i} G_s> - C\int_0^T ds \int_{D^d} du G(s,u)^2.$$
Note that
\begin{equation}\label{energyineq}
\sup_{G\in C_c^\infty(\Omega_T)} \{Q_{i,C}^G\} = \frac{\mathcal{ Q}_i(\pi)}{4C}.
\end{equation}
\begin{lemma}\label{lemaenergy}
There exists a constant $C_0=C_0(\kappa)>0$, such that for every $i=1,\ldots,d$, every $k=0,\ldots,d$, and every function $G$ in $C_c^\infty(\Omega_T)$
$$\limsup_{N\to\infty} \frac{1}{N^d} \log E_{\nu_\kappa^N} \left[ \exp\left\{ N^d \mathcal{Q}_{i,C_0}^G (\pi^{N,k})\right\}\right]\leq C_0.$$ 
\end{lemma}
\emph{Proof:} Applying Feynman-Kac's formula and using the same arguments in the proof of Lemma \ref{lemabound3}, we have that 
$$\frac{1}{N^d} \log E_{\nu_\kappa^N} \left[ \exp\left\{N \int_0^T ds \sum_{x\in D_N^{d}} (I_k(\eta_x(s))-I_k(\eta_{x-e_i}(s)))G(s,x/N) \right\} \right]$$
is bounded above by
$$\frac{1}{N^d} \int_0^T \lambda_s^N ds,$$
where $\lambda_s^N$ is equal to
$$\sup_f \Big\{\Big< N\!\!\!\!\!\sum_{x\in D_N^d}\!\!\!\!\! (I_k(\eta(x))-I_k(\eta(x-e_i)))G(s,x/N),f\Big>_{\!\!\!\nu_\kappa^N}\!\!\!\! + N^2<{\mathcal{L}}_N\sqrt{f},\sqrt{f}>_{\nu_\kappa^N}\Big\},$$
where the supremum is taken over all densities $f$ with respect to $\nu_\kappa^N$. By Proposition \ref{dirichlet}, the expression inside brackets is bounded above by
$$CN^d - \frac{N^2}{2} D_{\nu_\kappa^N}(f) + \sum_{x\in D_N^d}\left\{ NG(s,x/N)\int[I_k(\eta_x)-I_k(\eta_{x-e_i})]f(\eta)\nu_\kappa^N(d\eta) \right\}.$$
We now rewrite the term inside the brackets as
\begin{equation}\label{bracket}
\sum_{v\in\mathcal{V}} v_k \sum_{x\in D_N^d}\Big\{\int NG(s,x/N)  [\eta(x,v) - \eta(x-e_i,v)]f(\eta) \nu_\kappa^N(d\eta)\Big\}.
\end{equation}
After a simple computation, we may rewrite the terms inside the brackets of the above expression as
$$NG(s,x/N) \int [\eta(x,v) - \eta(x-e_i,v)] f(\eta) \nu_\kappa^N(d\eta)$$
\begin{eqnarray*}
&=&NG(s,x/N)\int \eta(x,v) f(\eta) \nu_\kappa^N(d\eta)\\
&-& NG(s,x/N) \int \eta(x,v) f(\eta^{x-e_i,x,v}) \frac{\gamma_{x-e_i,v}}{\gamma_{x,v}} \nu_\kappa^N(d\eta)\\
&=& NG(s,x/N) \int \eta(x,v)[f(\eta) - f(\eta^{x-e_i,x,v})] \nu_\kappa^N(d\eta)\\
&+& G\int \eta(x,v)f(\eta^{x-e_i,x,v})N\left[1-\frac{\gamma_{x-e_i,v}}{\gamma_{x,v}}\right]\\
&\leq& G(s,x/N)^2 \int f(\eta^{x-e_i,x,v})\nu_\kappa^N(d\eta)\\
&+& \frac{1}{4} \int\eta(x,v)f(\eta^{x-e_i,x,v}) \left[N\left(1-\frac{\gamma_{x-e_i},v}{\gamma_{x,v}}\right)\right]^2 \nu_\kappa^N(d\eta)\\
&+& N^2 \int \frac{1}{2} [\sqrt{f(\eta^{x-e_i,x,v})}- \sqrt{f(\eta)}]^2\nu_\kappa^N(d\eta)\\
&+& 2 G(s,x/N)^2 \int {\eta(x,v)}(\sqrt{f(\eta)}+\sqrt{f(\eta^{x-e_i,x,v})})^2\nu_\kappa^N(d\eta).
\end{eqnarray*}
Using the above estimate, we have that \eqref{bracket} is clearly bounded above by $C_1+C_1 G(s,x/N)^2$, by some positive constant $C_1 = C_1(\kappa)$, since $\gamma_{\cdot,v}$ is smooth and the fact that $f$ is a density with respect to $\nu_\kappa^N$. Thus, letting $C_0=C+C_1$, the statement of the Lemma holds. $\square$\\

It is well-known that $\mathcal{Q}(\pi)$ is finite if and only if $\pi$ has a generalized gradient, $\nabla\pi = (\partial_{u_1}\pi,\ldots,\partial_{u_d}\pi)$, and
$$\hat{\mathcal{Q}}(\pi) = \int_0^T\int_{D^d} du \|\nabla\pi_t(u)\|^2<\infty.$$
In which case, $\mathcal{Q}(\pi)=\hat{\mathcal{Q}}(\pi)$. Recall that the sequence $(Q_N)_N$ defined in the beginning of this section is tight. We have then the following proposition:
\begin{proposition}\label{sobolev}
Let ${Q^\ast}$ be any limit point of the sequence of measures $(Q^N)_N$. Then,
$$E_{Q^\ast} \left[ \int_0^T ds \left(\int_{D^d}\|\nabla \rho(s,u)\|^2 du\right)\right]<\infty,$$
and
$$E_{Q^\ast} \left[ \int_0^T ds \left(\int_{D^d}\|\nabla p_k(s,u)\|^2 du\right)\right]<\infty.$$
\end{proposition}
\emph{Proof:} We thus have to prove that the energy $\mathcal{Q}(\pi)$ is almost surely finite. Fix a constant $C_0>0$ satisfying the statement of Lemma \ref{lemaenergy}. Let $\{G_m:1\leq m\leq r\}$ be a sequence of functions in $C_0^\infty(\Omega_T)$ (the space of infinitely differentiable functions vanishing at the boundary) and $1\leq i\leq d$, and $0\leq k \leq d$, be integers. By the entropy inequality, there is a constant $C>0$ such that
$$E_{\mu^N}\left[\max_{1\leq m\leq r} \left\{\mathcal{Q}_{i,C_0}^{G_m} (\pi^{N,k})\right\}\right] \leq C + \frac{1}{N^d}\log E_{\nu_\kappa^N}\left[\exp\left\{ N^d \max_{1\leq m\leq r} \left\{ \mathcal{Q}_{i,C_0}^{G_m}(\pi^{N,k})\right\}\right\}\right].$$
Therefore, Lemma \ref{lemaenergy} together with the elementary inequalities 
$$\limsup_{N\to\infty} N^{-d}\log(a_N+b_N) \leq \limsup_{N\to\infty}\max\{\limsup_{N\to\infty} N^{-d} \log(a_N),\limsup_{N\to\infty} N^{-d}\log(b_N)\}$$ and $\exp\{\max\{x_1,\ldots,x_n\}\} \leq \exp(x_1)+\cdots+\exp(x_n)$ imply that
\begin{eqnarray*}
E_{Q^\ast} \left[\max_{1\leq m\leq r} \left\{\mathcal{Q}_{i,C_0}^{G_m}(\pi^{N,k})\right\}\right] &=&\lim_{N\to\infty}E_{\mu^N}\left[ \max_{1\leq m\leq r} \left\{\mathcal{Q}_{i,C_0}^{G_m}(\pi^{N,k})\right\}\right]\\
&\leq& C+C_0.
\end{eqnarray*}
Using this, the equation (\ref{energyineq}) and the monotone convergence theorem, we obtain the desired result. $\square$
\subsection{Proof of Theorem \ref{limhyd2}} Note that all limit points $Q^\ast$ of $(Q_N)_N$ are concentrated on absolutely continuous measures with respect to the Lebesgue measure since there is at most one particle per site, that is, 
$$Q^\ast\{\pi; \pi^k(du) = p_k(u)du,\hbox{~for all~}0\leq k\leq d\} = 1,$$
where $\pi^k$ denotes the $k$th component of $\pi$ and $p_0 = \rho$.

For $k=0,\ldots,d$, denote, again, by $\pi_t^{k,N}$ the empirical measure associated to the $k$th thermodynamic quantity:
\begin{equation*}
\pi_t^{k,N} = \frac{1}{N^d} \sum_{x \in D_N^d} I_k(\eta_x(t)) \delta_{x/N}.
\end{equation*}
Further, denote by $\pi_t^{k,N,b_1}$ and $\pi_t^{k,N,b_{N-1}}$ the empirical measures associated to the $k$th thermodynamic quantity
restricted to the boundaries:
$$
\pi_t^{k,N,b_i} = \frac{1}{N^{d-1}} \sum_{\substack{x \in D_N^d\\x_1=i}} I_k(\eta_x(t)) \delta_{x/N},
$$
for $i=1,N-1$. 

To compute ${\mathcal{ L}}_N <\pi_t^{k,N},H>$ for this process, we note that $\mathcal{ L}_N^c I_k(\eta_x)$ vanishes for $k=0,\ldots,d$, because the collision operator preserves local mass and momentum.

Since, in our definition of weak solution we considered test functions $H$ vanishing at the boundary, that is, $H(x) = 0$, if $x\in\{0,1\}\times\mathbb{T}^{d-1}$, we assume that $H$ vanishes at the boundary as well.

Now, we consider the martingale 
$$M_{t,k}^{N,H} = <\pi_t^{k,N},H>-<\pi_0^{k,N},H> - \int_0^t N^2\mathcal{ L}_N <\pi_s^{k,N},H>ds,$$
which can be decomposed into
\begin{eqnarray}
M_{t,k}^{N,H} &=& <\pi_t^{k,N},H> - <\pi_0^{k,N},H> - \int_0^t N^2\mathcal{ L}_N^{ex,1} <\pi_s^{k,N},H>ds \label{mart3}\\ 
&-& \int_0^t N^2\mathcal{ L}_N^{ex,2} <\pi_s^{k,N},H>ds - \int_0^t N^2\mathcal{ L}_N^{b} <\pi_s^{k,N},H>ds . \label{mart4}
\end{eqnarray}

We first prove that 
\begin{equation}\label{eqb1}
\int_0^t N^2 \mathcal{ L}_N^{b} <\pi_s^{k,N},H>ds
\end{equation}
vanishes as $N\to\infty$. A simple calculation shows that 
$$N^2\mathcal{ L}_N^{b} \eta(x,v) = N^2\left[\alpha_v(\tilde{x}/N)-\eta(x,v)\right],\;\;\hbox{if}\;\;x_1=1,$$
and
$$N^2\mathcal{ L}_N^{b} \eta(x,v) = N^2\left[\beta_v(\tilde{x}/N)-\eta(x,v)\right], \;\;\hbox{if}\;\;x_1=N-1.$$
Since $H$ vanishes on the boundary, $H((x+e_1)/N)=0$ if $x_1=N-1$, and $H((x-e_1)/N)=0$ if $x_1=0$. Then, we have the equalities $NH(x/N) = \partial_{x_1}^N H((x-e_1)/N)$, if $x_1=1$, and $NH(x/N) = -\partial_{x_1}^N H(x/N)$, if $x_1=N-1$. Therefore, we obtain
\begin{equation}\label{contaboundary}
\begin{array}{ccc}
N^2\mathcal{ L}_N^{b} <\pi^{k,N},H> \!\!\!&=&\!\!\! \frac{1}{N^{d-1}} \sum_{\substack{x\in D_N^d\\x_1=1}}\sum_{v\in{\mathcal{V}}} v_k[\alpha_v\left(\frac{\tilde{x}}{N}\right)-\eta(x,v)]\partial_{x_1}^N H\left(\frac{x-e_1}{N}\right)\\
&-& \frac{1}{N^{d-1}} \sum_{\substack{x\in D_N^d\\x_1=N-1}}\sum_{v\in{\mathcal{V}}} v_k[\beta\left(\frac{\tilde{x}}{N}\right)-\eta(x,v)]\partial_{x_1}^NH\left(\frac{x}{N}\right).
\end{array}
\end{equation}
We now use the last computation together with Lemma \ref{lemabound3} to conclude that \eqref{eqb1} vanishes as $N\to\infty$.

Further, after two summations by parts of the integrand on the right-hand term of (\ref{mart3}), we have that
\begin{eqnarray*}
\int_0^t N^2\mathcal{ L}_N^{ex,1} <\pi_s^{k,N},H>ds &=& \frac{1}{2}\int_0^t <\pi_s^{k,N},\Delta_N H> ds\\
&+&<\pi_t^{k,N,b_{N-1}},\partial_{u_1}^N H> - <\pi_t^{k,N,b_1},\partial_{u_1}^N H>, 
\end{eqnarray*}
and after one summation by parts on the right-hand term of (\ref{mart4}),
and noting again that $H$ vanishes at the boundaries, we have that
$$\int_0^t N^2\mathcal{ L}_N^{ex,2} <\pi_s^{k,N},H>ds =- \frac{1}{N^d}\int_0^t\sum_{j=1}^d \sum_{x\in\mathbb{ T}_N^d} (\partial_{u_j}^N H)\left(\frac{x}{N}\right) \tau_x W_{j,k}^{N,s} ds,$$
where $\tau_x$ stands for the translation by $x$ on the state space $X_N$ so that $(\tau_x \eta)(y,v)= \eta(x+y,v)$ for all $x,y\in\mathbb{ Z}^d$, $v\in\mathcal{ V}$, and $W_{j,k}^{N,s}$ is given by:
$$
W_{j,k}^{N,s} = \sum_{v\in\mathcal{ V}} v_k \sum_{z\in\mathbb{ Z}^d} p(z,v) z_j \eta_s(0,v)[1-\eta_s(z,v)],
$$
where $v_0=1$. Since $p(\cdot,v)$ is of finite range,
$$
E_{\mu_{\boldsymbol \lambda}^N} \left[ W_{j,k}^{N,s}\right] = \sum_{v\in\mathcal{ V}} v_k v_j \chi(\theta_v({\boldsymbol \lambda})),
$$
where $\chi(a) = a(1-a)$. 
Now, note that $E_{\nu_{\kappa}^N}(\eta(x,v)) = \alpha_v(x/N)$ if $x\in\{1\}\times\mathbb{T}_N^{d-1}$ and $E_{\nu_{\kappa}^N}(\eta(x,v)) = \beta_v(x/N)$
if $x\in\{N-1\}\times\mathbb{T}_N^{d-1}$.

We then apply Lemma \ref{replacement} to write the martingale in terms of the empirical measure. Further, we apply the replacement lemma for the boundary (Lemma \ref{lemabound3}) to obtain that all limit points satisfy the integral identity in the definition of weak solution of the problem \eqref{problema2}.  

Using the previous computations and the tightness of the sequence of measures $Q_N$ (for more details see \cite[Chapter 5]{KL}) we conclude that all limit points are concentrated on weak solutions of
\begin{equation*}
\partial_t (\rho,{\boldsymbol p}) + \sum_{v\in \mathcal{ V}} \tilde{v}\left[v\cdot \nabla\chi(\theta_v(\Lambda(\rho,{\boldsymbol p})))\right] = \frac{1}{2}\Delta (\rho,{\boldsymbol p}),
\end{equation*}
with boundary conditions, given in the trace sense, by
\begin{equation}\label{bound1}
(\rho,{\boldsymbol p})(t,x) = a(\tilde{x}), \hbox{~for~} x\in \{0\}\times\mathbb{T}^{d-1},
\end{equation}
and
\begin{equation}\label{bound2}
(\rho,{\boldsymbol p})(t,x) = b(\tilde{x}), \hbox{~for~} x\in \{1\}\times\mathbb{T}^{d-1},
\end{equation}
where $a(\cdot)$ and $b(\cdot)$ were defined in equation \eqref{funcaod}, and $v_0=1$. The uniqueness of weak solutions of the above equation implies that there is at most one limit point. Moreover, by Proposition \ref{sobolev}, each limit point of $(Q_N)_N$ is concentrated on a vector of measures with finite energy, that is: whose components have densities with respect to the Lebesgue measure that belong to the Sobolev space $H^1(D^d)$. This completes the proof of the theorem. $\square$

\section{Proof of the replacement lemma}\label{sec5}
As mentioned in the subsection \ref{subreplemma}, we only have to prove this result for the process without the boundary dynamics. In this case, we have a product invariant measure given by $\nu_{\rho,{\boldsymbol p}}^N$.

Let $\mu^N(T)$ be the Cesaro mean of $\mu^NS_t^N$, namely:
$$\mu^N(T) = \frac{1}{T} \int_0^T \mu^N S_t^N dt,$$
and let $\overline{f}_{T,k}^N$ be the Radon-Nikodym density of $\mu^N(T)$ with respect to $\nu_{\rho,{\boldsymbol p}}^N$. We have that the Dirichlet form of $\overline{f}_{T,k}^N$, $D_N(\overline{f}_{T,k}^N,\nu_{\rho,{\boldsymbol p}}^N)$, is bounded by $C N^{d-2}/2T$, where $C$ is some constant. Therefore, to prove the replacement lemma, it is enough to show that
$$\limsup_{\epsilon\to 0}\limsup_{N\to\infty} \sup_{D_N(f,\nu_{\rho,{\boldsymbol p}})<CN^{d-2}} \int \frac{1}{N^d} \sum_{x\in D_N^d} \tau_x V_{\epsilon N}^{j,k}(\eta(s)) f(\eta) \nu_{\rho,{\boldsymbol p}}^N(d\eta) = 0.$$

From now on we will simply write the Dirichlet form of a function $f$ with respect to the measure $\nu_{\rho,{\boldsymbol p}}^N$ as $D_N(f)$.

To prove the replacement lemma, we will prove the one and two block estimates:
\begin{lemma}
(One block estimate): For every constant $C > 0$, for $1\leq j\leq d$ and for $0\leq k\leq d$:
$$\limsup_{\ell\to\infty}\limsup_{N\to\infty}\sup_{D_N(f)\leq CN^{d-2}}\int \frac{1}{N^d} \sum_{x\in D_N^d} (\tau_x V_\ell^{j,k})(\eta) f(\eta) \nu_{\rho,{\boldsymbol p}}^N(d\eta) = 0,$$
where $V_\ell^{j,k}(\eta)$ was defined in Lemma \ref{replacement}.
\end{lemma}
\begin{lemma}\label{twoblocks}
(Two block estimate): For every constant $C > 0$, for $1\leq j\leq d$ and for $0\leq k\leq d$:
\begin{align*}
\limsup_{\ell\to\infty}\limsup_{\epsilon\to 0}\limsup_{N\to\infty} \sup_{D_N(f)\leq CN^{d-2}} \sup_{y\in\Lambda_{\epsilon N}} \int \frac{1}{N^d} \sum_{x\in D_N^d} \Big| {\boldsymbol I}^\ell (x+y)
- {\boldsymbol I}^{N\epsilon} (x)\Big| f(\eta) \nu_{\rho,{\boldsymbol p}}^N = 0,
\end{align*}
where ${\boldsymbol I}^\ell (x)$ was defined in subsection \ref{ensemble}.
\end{lemma}
\subsection{Proof of one block estimate}
We begin by noting that the exclusion rule and the fact that $\mathcal{ V}$ is finite prevents large densities or large momentum on ${\boldsymbol I}^\ell(0)$.

We have that the measure $\nu_{\rho,{\boldsymbol p}}^N$ is translation invariant. Therefore, we can write the sum on one block estimate as 
$$\int V_\ell^{j,k}(\eta) \left( \frac{1}{N^d}\sum_{x\in D_N^d} \tau_x f\right)(\eta) \nu_{\rho,{\boldsymbol p}}^N (d\eta) = \int V_\ell^{j,k}(\eta) \overline{f}(\eta) \nu_{\rho,{\boldsymbol p}}^N(d\eta),$$
where $\overline{f}$ stands for the space average of all translations of $f$:
$$\overline{f}(\eta) = \frac{1}{N^d} \sum_{x\in D_N^d} \tau_x f(\eta).$$

Denote by $X_\ell$ the configuration space $\left( \{0,1\}^\mathcal{ V}\right)^{\Lambda_\ell}$, by $\xi$ some configuration on $X_\ell$ and by $\nu_{\rho,{\boldsymbol p}}^\ell$ the product measure $\nu_{\rho,{\boldsymbol p}}^N$ restricted to $X_\ell$. For a density $f:X_N \to \mathbb{ R}_{+}$, $f_\ell$ stands for the conditional expectation of $f$ with respect to the $\sigma$-algebra generated by $\{ \eta(x,v): x\in\Lambda_\ell, v\in \mathcal{ V}\}$, that is obtained by integrating all coordinates outside this hypercube:
$$f_\ell(x_i) = \frac{1}{\nu_{\rho,{\boldsymbol p}}^\ell(\xi)} \int {\boldsymbol 1}_{\{\eta: \eta(z,v) = \xi(z,v),z\in\Lambda_\ell, v\in\mathcal{ V}\}} f(\eta) \nu_{\rho,{\boldsymbol p}}^N(d\eta),$$
for $\xi \in X_\ell$.

Since $V_\ell^{j,k}(\eta)$ depends on the configuration $\eta$ only through the occupation variables $\{ \eta(x,v): x\in\Lambda_\ell, v\in\mathcal{ V}\}$, in the last integral we can replace $\overline{f}$ by $\overline{f}_\ell$. In particular, to prove the lemma it is enough to show that
\begin{equation}\label{1bloco1}
\limsup_{\ell\to\infty}\limsup_{N\to\infty} \sup_{D_N(f)\leq CN^{d-2}} \int V_\ell^{j,k} (\xi) \overline{f}_\ell(\xi)\nu_{\rho,{\boldsymbol p}}^\ell (d\xi) = 0.
\end{equation}

We will now compute some estimates on the Dirichlet form. Let $<\cdot,\cdot>_\nu$ be the inner product in $L^2(\nu)$. For positive $f$, denote the Dirichlet form of f as:
\begin{eqnarray*}
D_N(f) &=& -<\sqrt{f}, (\mathcal{ L}_N^{ex}+\mathcal{ L}_N^c)f>_{\nu_{\rho,{\boldsymbol p}}^N}\\
       &=& -<\sqrt{f}, \mathcal{ L}_N^{ex,1}f>_{\nu_{\rho,{\boldsymbol p}}^N}-<\sqrt{f}, \mathcal{ L}_N^{ex,2}f>_{\nu_{\rho,{\boldsymbol p}}^N}-<\sqrt{f}, \mathcal{ L}_N^cf>_{\nu_{\rho,{\boldsymbol p}}^N}\\
       &:=& D_{N,1}(f) + D_{N,2}(f) + D_{N,c}(f).
\end{eqnarray*}
We have that
$$D_{N,1}(f) = \sum_{\substack{x,z\in D_N^d\\|x-z|=1}} I_{x,z}^{(1)}(f),$$
$$D_{N,2}(f) = \frac{1}{N}\sum_{x,z\in D_N^d}  I_{x,z}^{(2)}(f)$$
and
$$D_{N,c}(f) = \sum_{x\in D_N^d} I_{x}^{(c)}(f),$$
where
$$I_{x,z}^{(1)}(f) = \sum_{v\in\mathcal{ V}} \frac{1}{2} \int [\sqrt{f(\eta^{x,x+z,v})}-\sqrt{f(\eta)}]^2\nu_{\rho,{\boldsymbol p}}^N (d\eta),$$
$$I_{x,z}^{(2)}(f) = \sum_{v\in\mathcal{ V}} \int p(z,v) [\sqrt{f(\eta^{x,x+z,v})}-\sqrt{f(\eta)}]^2\nu_{\rho,{\boldsymbol p}}^N (d\eta)$$
and
$$I_{x,z}^{(c)}(f) = \sum_{q\in\mathcal{ Q}} \int p(x,q,\eta) [\sqrt{f(\eta^{x,q})}-f(\eta)]^2\nu_{\rho,{\boldsymbol p}}^N (d\eta).$$
Since the Dirichlet form is translation invariant and convex, we have that $D_N(\overline{f})\leq D_N(f).$

Now, let
$$D^\ell(h) = \sum_{\substack{x,z\in\Lambda_\ell\\|x-z|=1}} I_{x,z}^{\ell,(1)}(h)+\sum_{x,z\in\Lambda_\ell} \frac{1}{N} I_{x,z}^{\ell,(2)}(h) + \sum_{x\in\Lambda_\ell} I_{x}^{\ell,(c)}(h),$$
where each $I^{\ell,(i)}$ equals $I^{(i)}$ with $\nu_{\rho,{\boldsymbol p}}^N$ replacing $\nu_{\rho,{\boldsymbol p}}^\ell$. By using Schwarz inequality and the definition of $f_\ell$, we obtain that
$$I_{x,z}^{\ell,(1)}(\overline{f}_\ell) \leq I_{x,z}^{(1)}(\overline{f}), I_{x,z}^{\ell,(2)}(\overline{f}_\ell) \leq I_{x,z}^{(2)}(\overline{f}) \hbox{~and~} I_{x}^{\ell,(c)}(\overline{f}_\ell) \leq I_{x}^{(c)}(\overline{f})$$
for every $x,z\in\Lambda_\ell$. Therefore,
$$D^\ell(\overline{f}_\ell)\leq  \sum_{\substack{x,z\in\Lambda_\ell\\|x-z|=1}} I_{x,z}^{(1)}(\overline{f}_\ell)+\sum_{x,z\in\Lambda_\ell} \frac{1}{N} I_{x,z}^{(2)}(\overline{f}_\ell) + \sum_{x\in\Lambda_\ell} I_{x}^{(c)}(\overline{f}_\ell).$$
On the other hand, by translation invariance of $\overline{f}$, $I_{x,z}^{(1)}(\overline{f}) = I_{x+y,z+y}^{(1)}(\overline{f})$,
$I_{x,z}^{(2)}(\overline{f}) = I_{x+y,z+y}^{(2)}(\overline{f})$ and $I_{x}^{(c)}(\overline{f}) = I_{0}^{(c)}(\overline{f})$. Hence,
\begin{eqnarray*}
D^\ell(\overline{f}_\ell) &\leq& (2\ell+1)^d \sum_{i=1}^d I_{0,e_i}^{(1)}(\overline{f}) + \frac{(2\ell+1)^d}{N} \sum_{y\in\Lambda_\ell} I_{0,y}^{(2)} (\overline{f}) + (2\ell+1)^d I_0^{(c)}(\overline{f})\\
&\leq& \frac{(2\ell+1)^d}{N^d} (D_{N,1}(\overline{f}) + D_{N,2}(\overline{f}) + D_{N,c}(\overline{f})).
\end{eqnarray*}

Since the Dirichlet form is positive, $D_N(f)\leq CN^{d-2}$ implies that $D_{N,1}(f)\leq CN^{d-2},$ $D_{N,2}(f)\leq CN^{d-1}$ and $D_{N,c}(f)\leq CN^{d-2}$. Thus,
$$D^\ell (\overline{f}_\ell) \leq 3C(2\ell+1)^d N^{-2} := C_0(C,\ell)N^{-2}.$$
Therefore, the Dirichlet form of $\overline{f}_\ell$ vanishes as $N\uparrow\infty$. Hence, by \eqref{1bloco1}, to prove the one block estimate we must show that
\begin{equation}\label{onebl}
\limsup_{\ell\to\infty}\limsup_{N\to\infty} \sup_{D^\ell(f)\leq C_0(C,\ell)N^{-2}} \int V_\ell^{j,k}(\xi) f(\xi) \nu_{\rho,{\boldsymbol p}}^\ell (d\xi) = 0
\end{equation}
with the supremum carried over all densities with respect to $\nu_{\rho,{\boldsymbol p}}^N$.

We will now take the limit as $N\uparrow\infty$. To do so, we note that $V_\ell^{j,k}\leq C_1$, where $C_1$ is some constant, and therefore
$$\int_{X_\ell} V_\ell^{j,k}(\xi) f(\xi) \nu_{\rho,{\boldsymbol p}}^N(d\xi) \leq C_1.$$
This subset of $\mathcal{ M}_{+}(X_\ell)$ is compact for the weak topology, and since it is compact, for each $N$, there exists a density $f_N$ with Dirichlet form bounded by $C_0 N^{-2}$ that reaches the supremum. Let now $N_n$ be a subsequence such that 
$$\lim_{n\to\infty} \int V_\ell^{j,k} f_{N_n}(\xi) \nu_{\rho,{\boldsymbol p}}^\ell (d\xi) = \limsup_{N\to\infty} \int V_\ell^{j,k}(\xi) f_N(\xi) \nu_{\rho,{\boldsymbol p}}^\ell (d\xi).$$

To keep notation simple, assume, without loss of generality, that the sequences $N_n$ and $N$ coincide. By compactness, we can find a convergent subsequence $f_{N_n}$. Denote by $f_\infty$ the weak limit. Since the Dirichlet form is lower semicontinuous
$$D^\ell(f_\infty) = 0.$$
Moreover, by weak continuity,
$$\lim_{n\to\infty} \int V_\ell^{j,k} (\xi) f_{N_n}(\xi) \nu_{\rho,{\boldsymbol p}}^\ell (d\xi)= \int V_\ell^{j,k}(\xi)f_\infty(\xi) \nu_{\rho,{\boldsymbol p}}^\ell (d\xi).$$
In conclusion, expression (\ref{onebl}) is bounded above by
$$\limsup_{\ell\to\infty} \sup_{D^\ell(f) = 0} \int V_\ell^{j,k}(\xi) f(\xi) \nu_{\rho,{\boldsymbol p}}^\ell (d\xi).$$

We will now decompose along sets with a fixed number of conserved quantities. 

Recall that ${\mf V}_L$ is the set of all possible values of ${\boldsymbol I}^L(0)$ when $\eta$ runs over $(\{0,1\}^\mathcal{ V})^{\Lambda_L}$. Further, ${\mf V}_L$ is finite. Furthermore, consider for each ${\boldsymbol i}$ in ${\mf V}_L$ the canonical measure $\nu_{L,{\boldsymbol i}}$ defined in subsection \ref{ensemble}; and moreover, recall that
$$\nu_{\Lambda_L,{\boldsymbol i}} (\cdot ) = \mu^{\Lambda_L}_{{\boldsymbol \lambda}} \left(  \cdot \left\vert {\boldsymbol I}^L ={\boldsymbol i}\right.\right).$$

A probability density with Dirichlet form equal to zero is constant on each set with a fixed number of conserved quantities. It is convenient therefore to decompose each density $f$ along these sets. Thus
$$\int V_\ell^{j,k}(\xi) f(\xi) \nu_{\rho,{\boldsymbol p}}^N (d\xi) = \sum_{{\boldsymbol j}\in{\mf V}_\ell} T_{\boldsymbol j}(f) \int V_\ell^{j,k} \nu_{\ell,{\boldsymbol j}}(d\xi),$$
where,
$$T_j(f) = \int {\boldsymbol 1}_{\mathcal{ H}_\ell({\boldsymbol j})} f(\xi) \nu_{\rho,{\boldsymbol p}}^\ell(d\xi).$$

Since $\sum_{{\boldsymbol j}\in\mathcal{ H}_\ell({\boldsymbol j})} T_j(f) = 1$, to conclude the proof of the one block estimate, we must show that
$$\limsup_{\ell\to\infty} \sup_{{\boldsymbol j}\in{\mf V}_\ell} \int V_\ell^{j,k}(\xi) \nu_{\ell,{\boldsymbol j}}(d\xi) =0.$$
Since the measure $\nu_{\ell,{\boldsymbol j}}$ is concentrated on configurations with conserved quantity ${\boldsymbol j}$, the last integral equals
$$\int \left| \frac{1}{(2\ell+1)^d}\sum_{y\in\Lambda_\ell} \sum_{v\in\mathcal{ V}} v_k \sum_z p(z,v) z_j \tau_y(h(\xi,z,v))\!\!-\!\!\sum_{v\in\mathcal{ V}} v_j v_k E_{\nu_{\boldsymbol j}^\ell} [h(\xi,e_1,v)]\right|\nu_{\ell,{\boldsymbol j}}(d\xi),$$
where $h(\xi,z,v) = \xi(0,v)(1-\xi(z,v))$.

Fix some positive integer $n$, that shall increase to infinity after $\ell$. Decompose the set $\Lambda_\ell$ in cubes of length $2k+1$. Consider the set $A = \left\{ (2n+1)x, x\in\mathbb{ Z}^d\right\}\cap\Lambda_{\ell-n}$ and enumerate its elements: $A = \{x_1,\ldots,x_q\}$ in such a way that $|x_i|\leq |x_j|$ for $i\leq j$. For $1\leq i \leq q$, let $B_i = x_i + \Lambda_n$. Note that $B_i\cap B_j=\emptyset$ if $i\neq j$ and that $\cup_{1\leq i\leq q} B_i \subset \Lambda_\ell$. Let $B_0 = \Lambda_\ell - \cup_{1\leq i\leq q} B_i$. By construction $|B_0| \leq Kn\ell^{d-1}$ for some universal constant $K$. The previous integral is bounded above by
$$\sum_{i=0}^q \frac{|B_i|}{|\Lambda_\ell|}\!\int\! \left|\sum_{v\in\mathcal{ V}} v_k\!\!\left(\!\!\frac{1}{|B_i|}\sum_{y\in B_i} \sum_z p(z,v) z_j \tau_y(h(\xi,z,v))- v_j E_{\nu_{\boldsymbol j}^\ell} [h(\xi,e_1,v)]\right)\!\!\right|\nu_{\ell,{\boldsymbol j}}(d\xi).$$
Since $|B_0|\leq Kn\ell^{d-1}$, $\sum_v v_k \xi(0,v)(1-\xi(z,v))$ has mean $\sum_v v_k \chi(\theta_v(\Lambda({\boldsymbol j})))$, and $\left|\sum_{z\in B_i} p(z,v)z_j\right|$ is bounded, the sum is equal to
$$\frac{|\Lambda_n|}{|\Lambda_\ell|}\!\!\sum_{i=0}^q \!\int\! \left|\!\sum_{v\in\mathcal{ V}}\! v_k\left( \!\frac{1}{|B_n|}\sum_{y\in B_i} \sum_z p(z,v) z_j \tau_y(h(\xi,z,v))- v_j E_{\nu_{\boldsymbol j}^\ell} [h(\xi,e_1,v)]\!\right)\!\!\right|\nu_{\ell,{\boldsymbol j}}(d\xi)$$
plus a term of order $\mathcal{ O}(n/\ell)$. Since the distribution of $\{\xi(z,v); z\in B_i, v\in \mathcal{ V}\}$ does not depend on $i$, the previous sum is equal to
$$\int \left|\sum_{v\in\mathcal{ V}} v_k\left( \frac{1}{(2n+1)^d}\sum_{y\in \Lambda_n} \sum_z p(z,v) z_j \tau_y(h(\xi,z,v))- v_j E_{\nu_{\boldsymbol j}^\ell} [h(\xi,e_1,v)]\right)\right|\nu_{\ell,{\boldsymbol j}}(d\xi)$$
plus a term of order $\mathcal{ O}(n/\ell)$.

Now, let $\mu_{\boldsymbol \lambda}$ be the product measure on $\left( \{0,1\}^\mathcal{ V}\right)^{\mathbb{ Z}^d}$ with marginals given by 
$$\mu_{\boldsymbol \lambda} \{\eta: \eta(x,\cdot) = \xi\} = m_{\boldsymbol \lambda}(\xi),$$
for each $\xi\in \{0,1\}^\mathcal{ V}$ and $x\in \mathbb{ Z}^d$. Therefore, $E_{\nu_{\boldsymbol j}^\ell} [\xi(0,v)(1-\xi(e_1,v))] = E_{\nu_{\boldsymbol j}} [\xi(0,v)(1-\xi(e_1,v))]$, where $\nu_{\boldsymbol j} = \mu_{\Lambda({\boldsymbol j})}$. Moreover, if in the equivalence of ensembles we choose $L = L(\ell) = \lfloor C(\ell,\mathcal{ V})\rfloor$, where $C(\ell,\mathcal{ V})$ is the constant given in the equivalence of ensembles, we can replace the canonical measure by the grand canonical measure paying a price of order $o_\ell(1)$. Therefore, we can write the previous integral as
$$\int \left|\sum_{v\in\mathcal{ V}} v_k\left( \frac{1}{(2n+1)^d}\sum_{y\in \Lambda_n} \sum_z p(z,v) z_j \tau_y(h(\xi,z,v))- v_j E_{\nu_{\boldsymbol j}} [h(\xi,e_1,v)]\right)\right|\nu_{\boldsymbol j}^\ell(d\xi)$$
plus a term of order $o_\ell(1)$. We now note that $\nu_{\boldsymbol j}$ equals $\nu_{\boldsymbol j}^\ell$ on $\Lambda_\ell$. Then, the integral can be written as
$$\int \left|\sum_{v\in\mathcal{ V}} v_k\left( \frac{1}{(2n+1)^d}\sum_{y\in \Lambda_n} \sum_z p(z,v) z_j \tau_y(h(\xi,z,v))- v_j E_{\nu_{\boldsymbol j}} [h(\xi,e_1,v)]\right)\right|\nu_{\boldsymbol j}(d\xi)$$
plus a term of order $o_\ell(1)$. Let now,
$$g_{\boldsymbol j}(\xi) = \left|\sum_{v\in\mathcal{ V}} v_k\left( \frac{1}{(2n+1)^d}\sum_{y\in \Lambda_n} \sum_z p(z,v) z_j \tau_y(h(\xi,z,v))- v_j E_{\nu_{\boldsymbol j}} [h(\xi,e_1,v)]\right)\right|,$$
but we know that $E_{\nu_{\boldsymbol j}}[h(\xi,e_1,v)] = \chi(\theta_v(\Lambda(\bs j)))$, then,
$$g_{\boldsymbol j}(\xi) = \left|\sum_{v\in\mathcal{ V}} v_k\left( \frac{1}{(2n+1)^d}\sum_{y\in \Lambda_n} \sum_z p(z,v) z_j \tau_y(h(\xi,z,v))- v_j \chi(\theta_v(\Lambda(\bs j)))\right)\right|.$$
Now, $\left(\{0,1\}^\mathcal{ V}\right)^{\mathbb{ Z}^d}$ is compact on the product topology, and also, all the marginals of $\nu_{\boldsymbol j}$ converge to the marginals of $\nu_{\rho,{\boldsymbol p}}$, when ${\boldsymbol j}\to(\rho,{\boldsymbol p})$ as $\ell\to\infty$. Then, $\nu_{\boldsymbol j}$ converges weakly to $\nu_{\rho,{\boldsymbol p}}$. Further, since $g_{\boldsymbol j}(\xi) \to g_{\rho,{\boldsymbol p}}(\xi)$ for every $\xi$, we have from Theorem 5.5 of Billingsley \cite{B}, that
$$\int g_{\boldsymbol j}(\xi) \nu_{\boldsymbol j}(d\xi) \stackrel{\ell\to\infty}{\longrightarrow} \int g_{\rho,{\boldsymbol p}}(\xi) \nu_{\rho,{\boldsymbol p}}(d\xi),$$
this convergence being uniform on compact subsets of $\mathbb{ R}_{+}\times \mathbb{ R}^d$. Then, since the remainder term is $o_\ell(1)$, the limit as $\ell\to\infty$ and ${\boldsymbol j}\to(\rho,{\boldsymbol p})$ is
$$\int \left| \frac{1}{(2n+1)^d}\!\!\sum_{y\in\Lambda_n} \sum_{v\in\mathcal{ V}} v_k \sum_z z_j p(z,v)\tau_y(h(\xi,z,v))\!\! -\!\! \sum_{v\in\mathcal{ V}} v_j v_k \chi(\theta_v(\Lambda(\rho,{\boldsymbol p})))\right| \nu_{\rho,{\boldsymbol p}}(d\xi).$$
On the other hand, as $k\uparrow\infty$, by the law of large numbers, this integral converges to $0$. \\
Therefore, the one block estimate is proved. $\square$
\subsection{Proof of the two block estimate}
To prove the two block estimate, it is enough to show that
\begin{align}
\limsup_{\ell\to\infty}\limsup_{\epsilon\to 0}\limsup_{N\to\infty} \sup_{D_N(f)\leq CN^{d-2}} &\sup_{y\in (\Lambda_{\epsilon N}\setminus \Lambda_\ell)} \int \frac{1}{N^d} \sum_{x\in D_N^d} \Big| {\boldsymbol I}^\ell (x)\nonumber\\
&- {\boldsymbol I}^{\ell} (x+y)\Big| f(\eta)\nu_{\rho,{\boldsymbol p}}^N(d\eta) = 0.\label{int1twobl}
\end{align}
As for the one block estimate, we can rewrite this integral as 
$$\int \left| {\boldsymbol I}^\ell (0) - {\boldsymbol I}^{\ell} (y)\right| \overline{f}(\eta)\nu_{\rho,{\boldsymbol p}}^N(d\eta),$$
where $\overline{f}$ stands for the average of all space translations of $f$. ${\boldsymbol I}^\ell (0)$ and ${\boldsymbol I}^\ell (y)$ depend of the configuration $\eta$ only through the occupation variables $\{\eta(x,v): x\in\Lambda_{y,\ell}, v\in\mathcal{ V}\}$, where
$$\Lambda_{y,\ell} = \{-\ell,\ldots,\ell\}^d\cup [y+\{-\ell,\ldots,\ell\}^d].$$
We now introduce some notation. For positive integer $\ell$, let $X^{2,\ell}$ denote the configuration space $\left(\{0,1\}^\mathcal{ V}\right)^{\Lambda_\ell}\times \left(\{0,1\}^\mathcal{ V}\right)^{\Lambda_\ell}$, $\xi = (\xi_1,\xi_2)$ the configurations of $X^{2,\ell}$ and the product measure $\nu_{\rho,{\boldsymbol p}}^N$ restricted to $X^{2,\ell}$ (which does not depend on $N$) by $\nu_{\rho,{\boldsymbol p}}^{2,\ell}$. Denote by $f_{y,\ell}$ the conditional expectation of $f$ with respect to the $\sigma$-algebra generated by $\{\eta(x,v): x\in \Lambda_{y,\ell}, v\in\mathcal{ V}\}.$

Since ${\boldsymbol I}^\ell (0)$ and ${\boldsymbol I}^\ell (y)$ depend on $\eta(x,v)$, for $x\in\Lambda_{y,\ell}$ and $v\in\mathcal{ V}$, we may replace $\overline{f}$ by $\overline{f}_{y,\ell}$, and then, we can rewrite (\ref{int1twobl}) as
\begin{align*}
\limsup_{\ell\to\infty}\limsup_{\epsilon\to 0}\limsup_{N\to\infty} \sup_{D_N(f)\leq CN^{d-2}} \sup_{y\in (\Lambda_{\epsilon N}\setminus \Lambda_\ell)} \int \frac{1}{N^d} \sum_{x\in D_N^d} \Big| {\boldsymbol E}_1^\ell (0)
- {\boldsymbol E}_2^{\ell} (0)\Big| \overline{f}_{y,\ell}(\xi)\nu_{\rho,{\boldsymbol p}}^{2,\ell}(d\xi) = 0,
\end{align*}
where
$${\boldsymbol E}_i^{\ell}(x) = \frac{1}{|\Lambda_\ell|} \sum_{z\in x+\Lambda_\ell} {\boldsymbol I}(\xi_{iz}).$$
Now, we need to obtain information concerning the density $\overline{f}_{y,\ell}$ from the bound on the Dirichlet form of $f$. Then, let $D^{2,\ell}$ be the Dirichlet form defined on positive densities $h:X^{2,\ell}\to\mathbb{R}_{+}$ by
$$D^{2,\ell}(h) = I_{0,0}^\ell (h) + D_1^{\ell}(h) + D_2^\ell (h),$$
where,
\begin{eqnarray*}
D_1^\ell(h)\!\!\!\!\!\! &=&\!\!\!\!\!\! \sum_{v\in\mathcal{ V}} \int \left[ \sum_{\substack{x,z\in\Lambda_\ell\\|x-z|=1}}\frac{1}{2}+\frac{1}{N}\sum_{x,z\in\Lambda_\ell} p(z,v)\right] \left[\sqrt{h(\xi_1^{x,x+z,v},\xi_2)}-\sqrt{h(\xi)}\right]^2\nu_{\rho,{\boldsymbol p}}^{2,\ell}(d\xi)\\
&+&\sum_{x\in \Lambda_\ell}\sum_{v\in\mathcal{ V}}\int p(x,q,\xi_1)\left[\sqrt{h(\xi_1^{x,q},\xi_2)} - \sqrt{h(\xi)}\right]^2\nu_{\rho,{\boldsymbol p}}^{2,\ell}(d\xi),
\end{eqnarray*}
\begin{eqnarray*}
D_2^\ell(h)\!\!\!\!\!\! &=&\!\!\!\!\!\! \sum_{v\in\mathcal{ V}} \int \left[ \sum_{\substack{x,z\in\Lambda_\ell\\|x-z|=1}}\frac{1}{2}+\frac{1}{N}\sum_{x,z\in\Lambda_\ell} p(z,v)\right] \left[\sqrt{h(\xi_1,\xi_2^{x,x+z,v})}-\sqrt{h(\xi)}\right]^2\nu_{\rho,{\boldsymbol p}}^{2,\ell}(d\xi)\\
&+&\sum_{x\in \Lambda_\ell}\sum_{v\in\mathcal{ V}}\int p(x,q,\xi_1)\left[\sqrt{h(\xi_1,\xi_2^{x,q})} - \sqrt{h(\xi)}\right]^2\nu_{\rho,{\boldsymbol p}}^{2,\ell}(d\xi),
\end{eqnarray*}
and,
\begin{eqnarray*}
I_{0,0}^\ell(h) &=& \sum_{v\in\mathcal{ V}} \int \left[ \sum_{|z|=1}\frac{1}{2}+\frac{1}{N}p(z,v)\right] \left[\sqrt{h(\xi_1^{0,-,v},\xi_2^{0,+,v})}-\sqrt{h(\xi)}\right]^2\nu_{\rho,{\boldsymbol p}}^{2,\ell}(d\xi)\\
&+&\sum_{v\in\mathcal{ V}}\int p(0,q,\xi_1)\left[\sqrt{h(\xi_1^{0,q},\xi_2)} - \sqrt{h(\xi)}\right]^2\nu_{\rho,{\boldsymbol p}}^{2,\ell}(d\xi)\\
&+& \sum_{v\in\mathcal{ V}} \int \left[ \sum_{|z|=1}\frac{1}{2}+\frac{1}{N}p(z,v)\right] \left[\sqrt{h(\xi_1^{0,+,v},\xi_2^{0,-,v})}-\sqrt{h(\xi)}\right]^2\nu_{\rho,{\boldsymbol p}}^{2,\ell}(d\xi)\\
&+&\sum_{v\in\mathcal{ V}}\int p(0,q,\xi_2)\left[\sqrt{h(\xi_1,\xi_2^{0,q})} - \sqrt{h(\xi)}\right]^2\nu_{\rho,{\boldsymbol p}}^{2,\ell}(d\xi),
\end{eqnarray*}
where
$$\xi_i^{0,\pm,v}(x,w) = \left\{\begin{array}{cc}
\xi_i(0,v)\pm 1,& \hbox{if~}x=0\hbox{~and~}w=v,\\
\xi_i(x,w),&\hbox{otherwise}.
\end{array}\right.$$

This Dirichlet form corresponds to an interacting particle system on $(\mathcal{ V}\times \Lambda_\ell)\times(\mathcal{ V}\times \Lambda_\ell)$, where particles evolve according to an exclusion process with collisions among velocities on each coordinate and where particles from the origin of one of the coordinates at some velocity can jump to the origin of the other at this velocity and vice-versa.

Using the same idea as for the one-block estimate, we can prove that 
$$D_1^\ell(\overline{f}_{y,\ell})\leq D_N(\overline{f})\hbox{~~and~~}D_2^\ell(\overline{f}_{y,\ell})\leq D_N(\overline{f}),$$
and hence,
$$D_1^\ell(\overline{f}_{y,\ell}) + D_2^{\ell}(\overline{f}_{y,\ell})\leq 2C_0N^{-2},$$
for every density $f$ with Dirichlet form $D_N(f)$ bounded by $CN^{d-2}$. It remains to be shown that we can also estimate the Dirichlet form $I_{0,0}^\ell(\overline{f}_{y,\ell})$ by the Dirichlet form of $f$.

We begin by noting that 
$$I_{0,0}^\ell(h) = I_{0,0}^{\ell,1}(h)+I_{0,0}^{\ell,2}(h),$$
where,
\begin{align*}
I_{0,0}^{\ell,1}(h) = \sum_{v\in\mathcal{ V}}\left[\sum_{|z|=1}\frac{1}{2}+\frac{1}{N}p(z,v)\right]\Bigg[&\int\left[\sqrt{h(\xi_1^{0,-,v},\xi_2^{0,+,v})} - \sqrt{h(\xi)}\right]^2\\
&+ \left[\sqrt{h(\xi_1^{0,+,v},\xi_2^{0,-,v})} - \sqrt{h(\xi)}\right]^2\nu_{\rho,{\boldsymbol p}}^{2,\ell}(d\xi)\Bigg],
\end{align*}
and
\begin{eqnarray*}
I_{0,0}^{\ell,2}(h) &=& \sum_{v\in\mathcal{ V}} \int p(0,q,\xi_1)\left[\sqrt{h(\xi_1^{0,q},\xi_2)} - \sqrt{h(\xi)}\right]^2\nu_{\rho,{\boldsymbol p}}^{2,\ell}(d\xi)\\
&+& \sum_{v\in\mathcal{ V}}\int p(0,q,\xi_2)\left[\sqrt{h(\xi_1,\xi_2^{0,q})} - \sqrt{h(\xi)}\right]^2\nu_{\rho,{\boldsymbol p}}^{2,\ell}(d\xi).
\end{eqnarray*}
Then, a simple calculation shows that
$$I_{0,0}^{\ell,2}(\overline{f}_{y,\ell}) \leq 2I_0^{(c)}(\overline{f}),$$
and therefore $I_{0,0}^{\ell,2}(\overline{f}_{y,\ell})$ is also of order $N^{-2}$. We then have to obtain a bound for $I_{0,0}^{\ell,1}(\overline{f}_{y,\ell})$. 

Following the same lines used to prove that $I_{x,z}^{\ell,(j)}(\overline{f}_\ell) \leq I_{x,z}^{(j)}(\overline{f})$ in the proof of the one block estimate, for $j=1,2,c$, we have that each density $f$, with respect to $\nu_{\rho, {\boldsymbol p}}^N$, $I_{0,0}^{\ell,1}(\overline{f}_{y,\ell})$, is bounded above by:

\begin{equation}\label{int2twobl}
2 \sum_{v\in\mathcal{ V}}\left[\sum_{|z|=1}\frac{1}{2}+\frac{1}{N}p(z,v)\right] \int \left[\sqrt{\overline{f}(\eta^{0,y,v})}-\sqrt{\overline{f}(\eta)}\right]^2\nu_{\rho,{\boldsymbol p}}^N(d\eta).
\end{equation}

Let $(x_k)_{0\leq k\leq\||y\||}$ be a path from the origin to $y$, that is, a sequence of sites such that the first one is the origin, the last one is $y$ and the distance between two consecutive sites is equal to 1:
$$
x_0 = 0, x_{\||y\||}=y\hbox{~and~} |x_{k+1}-x_k|=1\hbox{~for every~}0\leq k\leq |||y|||-1,
$$
$|||\cdot|||$ is the sum norm:
$$
|||(y_1,\ldots,y_d)||| = \sum_{1\leq i\leq d}|y_i|.
$$
Let $\tau_{x_1}\cdots\tau_{x_i}(\eta)$ be the sequence of nearest neighbor exchanges that represents the path along $x_1,\ldots,x_i$. Then,
by using the telescopic sum 
$$
\sqrt{f(\eta^{0,y,v})} - \sqrt{f(\eta)} = \sum_{k=0}^{|||y|||-1}\left(\sqrt{f(\prod_{i=1}^k \tau_{x_i}(\eta))}-
\sqrt{f(\prod_{i=1}^{k-1} \tau_{x_i}(\eta))}\right)
$$
and the Cauchy-Schwarz inequality
$$
\left(\sum_{k=0}^{|||y|||-1} a_k\right)^2 \leq |||y||| \sum_{k=0}^{|||y|||-1} a_k^2,
$$
we obtain that (\ref{int2twobl}) is bounded by
$$
2\sum_{v\in\mathcal{ V}} \left[\sum_{|z|=1}\frac{1}{2}+\frac{1}{N}p(z,v)\right]\!\!|||y|||\!\!\sum_{k=0}^{|||y|||-1}\!\!\!\left[\sqrt{\overline{f}(\prod_{i=1}^k \tau_{x_i}(\eta))} - \sqrt{\overline{f}(\prod_{i=1}^{k-1} \tau_{x_i}(\eta))}\right]^2\!\! \nu_{\rho,{\boldsymbol p}}^N(d\eta)
$$
$$\leq 2\cdot 2\cdot 2^d |||y||| \sum_{k=0}^{|||y|||-1}I_{x_k,x_{k+1}}^{(1)}(\overline{f}).$$
Since $\overline{f}$ is translation invariant, for each $k$, $I_{x_k,x_{k+1}}^{(1)}(\overline{f}) = I_{x_k+z,x_{k+1}+z}^{(1)}(\overline{f})$
for all $z\in{\mathbb{Z}}^d$. Hence, $I_{x_k,x_{k+1}}^{(1)}(\overline{f})\leq N^{-d} D_N(f)$. In particular,
$$I_{0,0}^{\ell,1}(\overline{f}_{y,\ell}) \leq 2^{d+2}|||y|||^2 N^{-d}D_N(f).$$
Recall that $y\in\Lambda_{\epsilon N}$, and hence $|y|\leq 2N\epsilon$, $|\cdot |$ is the max norm. Then, $|||y|||\leq d|y|\leq 2dN\epsilon.$ Since the Dirichlet form is assumed to be bounded by $CN^{d-2}$, we have proved that 
$$I_{0,0}^{\ell,1}(\overline{f}_{y,\ell})\leq 2^{d+4}d^2C\epsilon^2.$$
We have, therefore, proved that for every density $f$ with Dirichlet form bounded by $CN^{d-2}$ and for every $d$-dimensional integer with max norm between $2\ell$ and $2N\epsilon$,
$$D^{2,\ell}(\overline{f}_{y,\ell}) \leq C_2(C,d,\ell)\epsilon^2.$$
We can now restrict ourselves to densities $f$ such that $D^{2,\ell}(\overline{f}_{y,\ell})\leq C_2\epsilon^2$, that vanishes as $\epsilon\downarrow 0$. In particular, to conclude the proof, it is enough to show that
$$\limsup_{\ell\to\infty}\limsup_{\epsilon\to 0}\sup_{D^{2,\ell}(f)\leq C_2\epsilon^2}\int |{\boldsymbol E}_1^{\ell}(0)-{\boldsymbol E}_2^\ell(0)| f(\xi)\nu_{\rho,{\boldsymbol p}}^{2,\ell}(d\xi) = 0,$$
this time, however, the supremum is taken over all densities with respect to $\nu_{\rho,{\boldsymbol p}}^{2,\ell}$. The rest of the proof follows the same lines as the ones in the one block estimate, beginning by decomposing the Dirichlet form along the sets having fixed conserved quantities and then applying the equivalence of ensembles.
Therefore, the two block estimate is proved. $\square$
\section{Uniqueness}\label{sec6}

To conclude the proof of the hydrodynamic limit, it remains to be proven the uniquenesses for the solutions of problems \eqref{problema2} and \eqref{problema2}. The strategy we used to prove this result was employed by Oleinik and Kruzhkov \cite{OK} and is due to Yu.A. Dubinskii.

Let $\nu$ and $\omega$ be two weak solutions to the problem (\ref{problema2}), 
corresponding to the same initial function $\nu_0$. Fix some $j=1,\ldots,d+1$, and let $H_j \in C^{1,2}\left([0,T]\times D^d\right)$ be such that $H_j(T,u)=0$, for all $u$. Then the integral identity for $\nu-\omega$ holds:
\begin{eqnarray}
\int_{0}^{T} dt\int_{D^d} du(\nu_j-\omega_j)\left[\partial_t H_j+\frac{1}{2} \sum_{1\leq i\leq d}\partial_{u_i}^2 H_j\right]\!\!
+\!\!\int_0^T dt \int_{D^d} du\! \sum_{v\in\mathcal{ V}} v_j(g_v(\nu)-g_v(\omega))\!\!\sum_{1\leq i \leq d} v_i \partial_{u_i}H_j =0,\label{uniqeq}
\end{eqnarray}
where $g_v(\nu) = \chi(\theta_v(\Lambda(\nu)))$, $\nu_j,\omega_j$ and $H_j$ are the components of $\nu,\omega$ and $H$, respectively. If $\nu_j=\omega_j$, we already have what we want, thus, suppose $\nu_j\neq\omega_j$. Introducing the
notation
$$\beta_v^j = \frac{g_v(\nu) - g_v(\omega)}{\nu_j-\omega_j},$$
we have that we can write (\ref{uniqeq}) as
\begin{equation}\label{uniqeq2}
\int_0^T dt \int_{D^d} du (\nu_j - \omega_j)\left[\partial_t H_j + \frac{1}{2} \sum_{1\leq i \leq d} \partial_{u_i}^2 H_j + \sum_{v \in \mathcal{ V}} v_j \beta_v^j \sum_{1\leq i\leq d} v_i \partial_{u_i} H_j\right] = 0.
\end{equation}
Now, let $\beta_v^{j,m}$ be a sequence of smooth functions which converge in $L^2([0,T]\times D^d)$ to $\beta_v^j$, as $m\to\infty$. We denote by $H^m_j(t,u)$ the classical solution of the equation
\begin{equation}\label{uniqeq3}
\partial_t H_j^m + \frac{1}{2} \sum_{1\leq i \leq d} \partial_{u_i}^2 H_j^m + \sum_{v \in \mathcal{ V}} v_j \beta_v^{j,m} \sum_{1\leq i\leq d} v_i \partial_{u_i} H_j^m = \Phi_j,
\end{equation}
$$H^m_j(T,u) = 0, H_j^m(0,u) = 0,$$
where $\Phi_j$ is a smooth function finite in $[0,T]\times D^d$. For more details on the solutions of partial differential equations of the parabolic type,
the reader is referred to Friedman \cite{F}, and for details on solutions of systems of linear partial differential equations of the parabolic type in general, the reader is referred to Lady\v{z}enskaja et al. \cite{L}.

Now, if we replace $H_j$ in (\ref{uniqeq2}) by $H^m_j$ and use (\ref{uniqeq3}), we obtain:
\begin{align}
\int_0^T &dt \int_{D^d} du (\nu_j - \omega_j)\Phi_j
+ \int_0^T dt \int_{D^d}du(\nu_j-\omega_j)\left[\sum_{v\in\mathcal{ V}} v_j (\beta_v^j- \beta_v^{j,m})\sum_{1\leq i \leq d} \partial_{u_i} H_j^m\right] = 0.\label{uniqeq4}
\end{align}

Finally, since we are in a compact domain and the coefficients $\beta_v^{j,m}$ are smooth, we have that there exists an $M>0$ such that $|H^m_j|\leq M$. Since these coeffiecients converge in $L^2([0,T]\times D^d)$, the constant $M$ may be taken to be independent of $m$. Multiplying \eqref{uniqeq3} by $H_j^m$, integrating over $[0,T]\times D^d$, and then integrating by parts, we have that
$$\int_0^Tdt\int_{D^d}du\sum_{i=1}^d \frac{1}{2}\left(\frac{\partial H_j^m}{\partial u_i}\right)^2 = \int_0^Tdt\int_{D^d}du \left( \sum_{v \in \mathcal{ V}} v_j \beta_v^{j,m} H_j^m \sum_{1\leq i\leq d} v_i \partial_{u_i} H_j^m - \Phi H_j^m\right)-\frac{1}{2}\int_{D^d}du (H_j^m)^2 .$$
On applying the elementary inequality $|ab|\leq \epsilon a^2 + b^2/(4\epsilon)$ and using that $|H_j^m|\leq M$, we obtain that
$$\int_0^Tdt\int_{D^d}du \sum_{i=1}^d\frac{1}{2}\left(\frac{\partial H_j^m}{\partial u_i}\right)^2\leq C,$$
where $C$ is a constant that may depend on $M$ and $\Phi$, but not on $m$. 

Therefore, by applying the Cauchy-Schwartz inequality and using that $\beta_v^{j,m}$ converges to $\beta_v^j$ in the $L^2$-norm, we see that the second term on the left-hand side of equation \eqref{uniqeq4} tends to zero as $m$ tends to infinity. This implies that for every $\varepsilon>0$ there exists $m$ such that the absolute value of the second term on the left-hand side of equation \eqref{uniqeq4} is less than $\varepsilon$. We, then, have obtained that
$$\forall\varepsilon>0: \left|\int_0^T dt \int_{\mathbb{ T}^d}du (\nu_j-\omega_j)\Phi_j\right|\leq \varepsilon,$$
and hence, for each $j=1,\ldots,d+1$, $\nu_j=\omega_j$. Therefore $\nu\equiv \omega$. $\square$

\section*{Acknowledgements} I would like to thank my PhD advisor, Claudio Landim, for suggesting this problem, for valuable
discussions and support. I also thank Jonathan Farfán for stimulating discussions on this topic, and Fábio Valentim for helping me
improve the overall quality of the paper.

\end{document}